\documentclass[twoside]{irmaems}
\usepackage[all]{xy}
\usepackage{amsmath}
\usepackage{amssymb}
\usepackage{latexsym}
\RequirePackage[T1]{fontenc}

\renewcommand\theequation{\arabic{section}.\arabic{equation}}

 \newtheorem{dfn}{Definition}[section] 
 \newtheorem{prop}[dfn]{Proposition}
 \newtheorem{cor}[dfn]{Corollary}
 \newtheorem{cjt}[dfn]{Conjecture}
 \newtheorem{thm}[dfn]{Theorem}
 \newtheorem{lem}[dfn]{Lemma}
 \newtheorem{ex}[dfn]{Example}     
 \newtheorem{exs}[dfn]{Examples}   
 \newtheorem{pb}[dfn]{Problem}

\newcommand{\rep}{\textrm{Rep}}
\newcommand{\obj}[1]{\textrm{Obj}(#1)}
\newcommand{\vect}{\textrm{Vect}_{\K}}
\newcommand{\End}[1]{\textrm{End}(#1)}

\newcommand{\Fun}[1]{\textrm{Fun}(#1)}
\newcommand{\Gr}[1]{\textrm{Gr}(#1)}
\newcommand{\Mor}[2]{\textrm{Hom}(#1,#2)}

\newcommand{\RR}{\mathbb{R}}
\newcommand{\Q}{\mathbb{Q}}
\newcommand{\N}{\mathbb{N}}
\newcommand{\Z}{\mathbb{Z}}
\newcommand{\C}{\mathbb{C}}
\newcommand{\K}{\textbf{k}}

\newcommand{\cat}{\mathcal C}
\newcommand{\dat}{\mathcal D}
\newcommand{\modu}{\mathcal M}

\newcommand{\al}{\mathfrak g}

\newcommand{\sld}{\mathfrak{sl}_2}

\newcommand{\aff}{\widehat{\mathfrak{sl}_2}}

\begin{document}

\title{Lectures on tensor categories}

\author{Damien Calaque and Pavel 
Etingof\thanks{The work of P.E.~was partially
supported by the NSF grant DMS-9988786. }}

\address{
IRMA (CNRS)\\
7 Rue Ren\'e Descartes, F-67084 Strasbourg, France\\
email:\,\tt{calaque@math.u-strasbg.fr}
\\[4pt]
\\
Department of Mathematics, MIT\\
Cambridge, MA 02138, USA\\
email:\,\tt{etingof@math.mit.edu}
}

\maketitle

\begin{abstract} We give a review of some recent developments in the theory 
of tensor categories. The topics include realizability of fusion rings, Ocneanu
rigidity, module categories, weak Hopf algebras, Morita theory for tensor 
categories, lifting theory, categorical dimensions, Frobenius-Perron 
dimensions classification of tensor categories.   \end{abstract}



\tableofcontents

\section*{Introduction}

This paper is based on the first author's lecture notes  of a series of four lectures given by the 
second author in June 2003 for the Quantum Groups seminar at the Institut de Recherche Mathématique 
Avancée in Strasbourg. It is about the structure and classification of tensor categories. 

We always work over an algebraically closed field $\K$. By a tensor category over $\K$, we mean 
an abelian rigid monoidal category in which the neutral object \textbf1 is simple (i.e., does not 
contain any proper subobject), the vector spaces $\Mor XY$ are finite dimensional and all 
objects are of finite length. 

The category of finite dimensional vector spaces $\vect$, the 
categories of finite dimensional representations of a group $G$, a Lie algebra $\al$, or a 
(quasi)Hopf algebra $H$ (respectively denoted $\rep G$, $\rep\al$ and $\rep H$), or the category 
of integrable modules over an affine Lie algebra $\hat\al$ with fusion product (which can also be 
obtained from quantum groups) are all tensor categories in this sense. 

Tensor categories appear in many areas of mathematics such as representation theory, quantum groups, 
conformal field theory (CFT) and logarithmic CFT, operator algebras, and topology (invariants of 
knots and 3-manifolds). The goal of this paper is to give an introduction to some recent developments 
in this subject. 

The paper is subdivided into four sections, each representing a single lecture. 

Section 1 introduces the main objects of the paper. We recall basic categorical definitions and 
results, fix the vocabulary, and give examples (for more details, we recommend the monographies 
\cite{BK,ML,K,T}). The end of the section is devoted to the problem of realizability of fusion 
rings: examples are given, and the Ocneanu rigidity conjecture is formulated. 

The goal of Section 2 is to prove the Ocneanu rigidity for fusion categories in characteristic zero. 
To do this, we introduce and discuss the notions of module categories and weak Hopf algebras. 
The more technical part of the proof is done at the end of the section. 

Section 3 is about three distinct subjects. We start with a closer look at module categories, 
discussing the notion of Morita equivalence for them, and applying general results to the 
representation theory of groups. Then we recall well-known facts about braided, ribbon and modular 
categories. Finally, the lifting theory is presented: it allows us to extend some results from 
characteristic zero to the positive characteristic case. 

Section 4 covers the theory of Frobenius-Perron dimension, and its applications to classification 
results for fusion categories. 

We end this paper with two interesting open problems. 

\medskip \noindent{\bf Remarks.} 
1. Being a set of lecture notes, this paper does not contain original results. 
Most of the results are taken from the papers \cite{ENO,EO1,O1,O2,O3} and references therein, 
including the standard texts on the theory of tensor categories. 

2. To keep this paper within bounds, we had to refrain from a thorough review of the history of the 
subject and of the original references, as well as from a detailed discussion of the preliminaries. 
We also often had to omit complete proofs. For all this material we refer the reader to books and 
papers listed in the bibliography. \\

\noindent\textbf{Acknowledgements. } The authors are grateful to the participants of the lectures 
-- P. Baumann, B. Enriquez, F. Fauvet, C. Grunspan, G. Halbout, C. Kassel, V. Turaev, and B. Vallette. 
Their interest and excitement made this paper possible. The second author is greatly indebted to 
D. Nikshych and V. Ostrik for teaching him much of the material given in these lectures. He is also 
grateful to IRMA (Strasbourg) for hospitality. His work was partially supported by the NSF grant 
DMS-9988786. 

\section{Finite tensor and fusion categories}

\subsection{Basic notions}

\subsubsection{Definitions}

Let $\mathcal C$ be a category. 

Recall that $\mathcal C$ is \emph{additive} over $\K$ if \\ 
(i) $\Mor XY$ is a (finite dimensional) $\K$-vector space for all $X,Y\in\obj{\mathcal C}$, \\ 
(ii) the map $\Mor YZ\times\Mor XY\to\Mor XZ,(\varphi,\psi)\mapsto\varphi\circ\psi$ is $\K$-linear 
for all $X,Y,Z\in\obj{\mathcal C}$, \\ 
(iii) there exists an object $\textbf0\in\obj{\mathcal C}$ such that 
$\Mor {\textbf0}X=\Mor X{\textbf0}=\textbf0$ for all $X\in\obj{\mathcal C}$, \\ 
(iv) finite direct sums exist. 

\medskip 
\noindent\textbf{Remark.} When we deal with functors between additive categories, we always 
assume they are also additive. 
\medskip 

Further, recall that an additive category $\cat$ is \emph{abelian} if \\
(i) every morphism $\phi: X\to Y$ has a kernel ${\rm Ker}\phi$ (an object $K$ together with a 
monomorphism $K\to X$) and a cokernel ${\rm Coker}\phi$ (an object $C$ together with an epimorphism 
$Y\to C$); \\
(ii) every morphism is the composition of an epimorphism followed by a monomorphism; \\
(iii) for every morphism $\varphi$ one has 
$\textrm{Ker}\varphi=0\implies\varphi=\textrm{Ker}(\textrm{Coker}\varphi)$ and 
$\textrm{Coker}\varphi=0\implies\varphi=\textrm{Coker}(\textrm{Ker}\varphi)$. 

It is known that $\cat$ is abelian if and only if 
it is equivalent to a full subcategory of the category of 
modules over a algebra. 
Recall also that $\mathcal C$ is \emph{monoidal} if there exists \\ 
(i) a bifunctor $\otimes:\mathcal C\times\mathcal C\to\mathcal C$, \\ 
(ii) a functorial isomorphism $\Phi:(-\otimes-)\otimes-\to-\otimes(-\otimes-)$, \\ 
(iii) an object \textbf1 (called the neutral object) and two functorial isomorphisms 
$$
\lambda:\textbf1\otimes-\to-\,,\,\,\,\mu:-\otimes\textbf1\to-
\quad(\textrm{the unit morphisms})\,, 
$$
such that for any two functors obtained from $-\otimes\cdots\otimes-$ by inserting \textbf1's and 
parentheses, all functorial isomorphisms between them composed of $\Phi^{\pm1}$'s, $\lambda^{\pm1}$'s 
and $\mu^{\pm1}$'s are equal. \\

\noindent\textbf{Remark. }In the spirit of the previous remark, for additive monoidal categories 
we assume that $\otimes$ is  biadditive. 

\begin{thm}[MacLane coherence, \cite{ML}] The data $(\mathcal C,\otimes,\Phi,\lambda,\mu)$ with (i), 
(ii) and (iii) is a monoidal category if and only if the following properties are satisfied: 
\begin{enumerate}
\item\textbf{Pentagon axiom. }The following diagram is commutative: \\
$\xymatrix{
((-\otimes-)\otimes-) \otimes-\ar[r]^{\Phi^{1,2,3}\otimes\textrm{id}}\ar[d]^{\Phi^{12,3,4}} 
& (-\otimes(-\otimes-))\otimes-\ar[r]^{\Phi^{1,23,4}} 
& -\otimes((-\otimes-)\otimes-)\ar[d]^{\textrm{id}\otimes\Phi^{2,3,4}} \\
(-\otimes-)\otimes(-\otimes-)\ar[rr]^{\Phi^{1,2,34}} & & -\otimes(-\otimes(-\otimes-))}$
\item\textbf{Triangle axiom. }The following diagram commutes: \\
$\xymatrix{
(-\otimes\textbf1) \otimes-\ar[r]^{\Phi_{-,\textbf1,-}}\ar[dr]_{\mu\otimes\textrm{id}} & 
-\otimes(\textbf1\otimes-)\ar[d]^{\textrm{id}\otimes\lambda} \\
& -\otimes-}$
\end{enumerate}
\end{thm}

A monoidal category is called {\it strict} if $(X\otimes Y)\otimes Z=X\otimes (Y\otimes Z)$, 
$\bold 1\otimes X=X\otimes \bold 1=X$, and the associativity and unit isomorphisms are equal to the 
identity. A theorem also due to Maclane (see \cite{K}) says that any monoidal category is equivalent 
to a strict one. In view of this theorem, we will always assume that the categories we are working 
with are strict, unless otherwise specified. 

Recall that a right dual for $X\in\obj{\mathcal C}$ is 
an object $X^*$ with two morphisms $e_X:X^*\otimes X\to\textbf1$ and $i_X:\textbf1\to X\otimes X^*$ 
(called the evaluation and coevaluation morphisms) satisfying the following two equations: \\
(i) $(\textrm{id}_X\otimes e_X)\circ(i_X\otimes \textrm{id}_X)=\textrm{id}_X$ and \\
(ii) $(e_X\otimes\textrm{id}_{X^*})\circ(\textrm{id}_{X^*}\otimes i_X)=\textrm{id}_{X^*}$. \\
A left dual $^*\!\!X$ with maps $e_X':X\otimes^*\!\!X\to\textbf1$ and 
$i_X':\textbf1\to^*\!\!X\otimes X$ is defined in the same way. \\
One can show that if it exists, the right (left) dual is unique up to a unique isomorphism compatible 
with evaluation and coevaluation maps. 

A monoidal category is called \emph{rigid} if any object has 
left and right duals. 

\begin{dfn}
\emph{A \emph{tensor category} is a rigid abelian monoidal category in which the object $\bold 1$ is 
simple and all objects have finite length. }
\end{dfn}
\begin{ex}\emph{
The category $\rep H$ of finite dimensional  representations of a quasi-Hopf algebra $H$ is a tensor 
category \cite{Dr}. This category is, in general, not strict (although it is equivalent to a strict 
one): its associativity isomorphism is given by the associator of $H$. 
}\end{ex}
\begin{prop}[\cite{BK}]
In a tensor category, the tensor product functor $\otimes$ is (bi)exact. 
\end{prop}

\subsubsection{The Grothendieck ring of a tensor category}

Let $\cat$ be a tensor category over $\K$. 
\begin{dfn}
\emph{The \emph{Grothendieck ring ${\rm Gr}(\cat)$} 
of $\cat$ is the ring whose basis over $\Z$ is the set 
of isomorphism classes of simple objects, with multiplication given by 
$$X\cdot Y=\sum_{Z~\textrm{simple}}N_{XY}^ZZ,$$
where $N_{XY}^Z=[X\otimes Y:Z]$ is the multiplicity (the number of occurences) of $Z$ in $X\otimes Y$ 
(which is well-defined by the Jordan-H\"older theorem). }
\end{dfn}

\begin{exs}\label{excat}
\emph{(i) $\cat=\rep_\C SL(2)$. Simple objects are highest weight representations $V_j$ 
(of highest weight $j\in \Bbb Z$), and the structure constants of the Grothendieck ring are given 
by the Clebsch-Gordan formula 
$$V_i\otimes V_j=\sum_{\substack{k= \vert i-j\vert \\ k\equiv i+j~\textrm{mod}~2}}^{i+j}V_k$$}

\emph{(ii) $\cat$ is the category of integrable modules (from category $\mathcal O$) over the affine 
algebra $\aff$ at level  $l$ with the fusion product 
$$V_i\otimes V_j=\sum_{\substack{k=\vert i-j\vert \\
k\equiv i+j~\textrm{mod}~2}}^{l-\vert i+j-l\vert}V_k$$
In this case $\Gr\cat$ is a Verlinde algebra. }

\emph{(iii) $\cat=\rep_\K\Fun G$ for a finite group $G$. Simple objects are evaluation modules $V_g$, 
$g\in G$, and $V_g\otimes V_h=V_{gh}$. So $\Gr\cat=\Z[G]$. \\
More generally, pick a 3-cocycle $\omega\in Z^3(G,\C^\times)$. To this cocycle we can attach a 
twisted version $\cat(G,\omega)$ of $\mathcal C$: all the structures are the same, except the 
associativity isomorphism which is given by $\Phi_{V_g,V_h,V_k}=\omega(g,h,k)\textrm{id}$ (and the 
morphisms $\lambda,\mu$ are modified to satisfy the triangle axiom). The cocycle condition 
$$\omega(h,k,l)\omega(g,hk,l)\omega(g,h,k)=\omega(gh,k,l)\omega(g,h,kl)$$
is equivalent to the pentagon axiom. Again, we have $\Gr{\cat(G,\omega)}=\Z[G]$. }

\emph{(iv) $\cat=\rep_\C S_3$. The basis elements (simple objects) are $\textbf1,\chi,V$, with 
product given by $\chi\otimes\chi=\textbf1$, $\chi\otimes V=V\otimes\chi=V$ and 
$V\otimes V=V\oplus\textbf1\oplus\chi$. }

\emph{(v) If $\cat=\rep G$ for $G$ a unipotent algebraic group over $\C$, then the unique simple 
object is $\bold 1$, hence $\Gr\cat=\Z$. In this case, the Grothendieck ring does not give a lot 
of information about the category because the category is not semisimple. }

\emph{(vi) $\cat=\rep H$ for the 4-dimensional Sweedler Hopf algebra $H$, which is generated by $g$ 
and $x$, with relations $gx=-xg$, $g^2=1$, $x^2=0$, and the coproduct $\Delta$ given by 
$\Delta g=g\otimes g$ and $\Delta x=x\otimes g+1\otimes x$. In this case the only simple objects are 
$\textbf1$ and $\chi$, with $\chi\otimes\chi=\textbf1$. }
\end{exs}

\subsubsection{Tensor functors}

Let $\cat$ and $\dat$ be two tensor categories. A functor $F:\cat\to\dat$ is called \emph{quasitensor} 
if it is exact and equipped with a functorial isomorphism $J:F(-\otimes-)\to F(-)\otimes F(-)$ and an 
isomorphism $u:F(\textbf1)\to\textbf1$. Such a functor defines a morphism of unital rings 
$\Gr\cat\to\Gr\dat$. \\
A quasitensor functor $F:\cat\to\dat$ is \emph{tensor} if the diagrams \\
$\xymatrix{
F((-\otimes-)\otimes-)\ar[r]^{J^{12,3}}\ar[d]^{F(\Phi_\cat)} & 
F(-\otimes-)\otimes F(-)\ar[r]^{J\otimes\textrm{id}} & 
(F(-)\otimes F(-))\otimes F(-)\ar[d]^{\Phi_\dat} \\
F(-\otimes(-\otimes-))\ar[r]^{J^{1,23}} & 
F(-)\otimes F(-\otimes-)\ar[r]^{\textrm{id}\otimes J} & 
F(-)\otimes(F(-)\otimes F(-))}$ \\
$\xymatrix{
F(\textbf1\otimes-)\ar[r]^{J_{\textbf1,-}}\ar[d]_{F(\lambda_\cat)} & 
F(\textbf1)\otimes F(-)\ar[d]^{u\otimes\textrm{id}} \\
F(-) & \textbf1\otimes F(-)\ar[l]^{\lambda_\dat}}$ and 
$\xymatrix{
F(-\otimes\textbf1)\ar[r]^{J_{-,\textbf1}}\ar[d]_{F(\mu_\cat)} & 
F(-)\otimes F(\textbf1)\ar[d]^{\textrm{id}\otimes u} \\
F(-) & F(-)\otimes\textbf1\ar[l]^{\mu_\dat}}$ \\
are commutative. \\
An \emph{equivalence of tensor categories} is a tensor functor which is also an equivalence of 
categories. 
\begin{ex}\emph{
Let $\omega,\omega'\in Z^3(G,\K^\times)$ and $\omega'/\omega= \textrm d\eta$ is a coboundary. Then 
$\eta$ defines a tensor structure on the identity functor $\cat(G,\omega') \to\cat(G,\omega)$: the 
coboundary condition 
$$\omega'(g,h,k) \eta(h,k)\eta(g,hk)=\eta(gh,k)\eta(g,h)\omega(g,h,k)$$
is equivalent to the commutativity of the previous diagram. Moreover, it is not difficult to see that 
this tensor functor is in fact an equivalence of tensor categories. Thus the fusion category 
$\cat(G,\omega)$, up to equivalence, depends only on the cohomology class of $\omega$. In particular, 
we may use the notation $\cat(G,\omega)$ when $\omega$ is not a cocycle but a cohomology class. 
}\end{ex}

\subsection{Finite tensor and fusion categories}

\subsubsection{Definitions and examples}

\begin{dfn}
\emph{An abelian category $\cat$ over $\K$ is said to be \emph{finite} if \\
(i) $\cat$ has finitely many (isomorphism classes of) simple objects, \\
(ii) any object has finite length, and \\
(iii) any simple object admits a projective cover. }
\end{dfn}
This is equivalent to the requirement that $\cat=\rep A$ as an abelian category for a finite 
dimensional $\K$-algebra $A$. 
\begin{dfn}
\emph{A \emph{fusion category} is a semisimple finite tensor category. }
\end{dfn}
\begin{exs}\emph{
In examples \ref{excat}, (i) is semisimple but not finite, (ii), (iii) and  (iv) are fusion, (v) is 
neither finite nor semisimple, and (vi) is finite but not semisimple. 
}\end{exs}

Recall that if $\cat$ and $\dat$ are two abelian categories over $\K$, then one can define their 
\emph{Deligne external product} $\cat\boxtimes\dat$. Namely, if $\cat=A$-Comod  and 
$\dat=B$-Comod are the categories of comodules over coalgebras $A$ and $B$ then 
$\cat\boxtimes\dat:=A\otimes B$-Comod. 

If $\cat$ and $\dat$ are semisimple, the Deligne product is simply the category whose simple objects 
are $X\boxtimes Y$ for simple $X\in\obj\cat$ and $Y\in\obj\dat$. If $\cat$ and $\dat$ are 
tensor/finite tensor/fusion categories then $\cat\boxtimes\dat$ also has a natural structure of a 
tensor/finite tensor/fusion category (in the semisimple case it is simply given by 
$(X\boxtimes Y)\otimes(X'\boxtimes Y'):=(X\otimes Y)\boxtimes(X'\otimes Y')$). 

\subsubsection{Reconstruction theory (Tannakian formalism)}

Let $H$ be a (quasi-)Hopf algebra and consider $\cat=\rep H$, the category of its finite dimensional 
representations. The forgetful functor $F:\cat\to\vect$ has a (quasi)tensor structure (the identity 
morphism). In addition, this functor is exact and faithful. A functor $\cat\to\vect$ with such 
properties ((quasi)tensor, exact, and faithful) is  called a {\it (quasi)fiber functor}. 

\emph{Reconstruction theory} tells us that every finite tensor category equipped with a (quasi)fiber 
functor is obtained in this way, i.e., can be realized as the category of finite dimensional 
representations of a finite dimensional (quasi-)Hopf algebra. 

Namely, let $(\cat,F)$ be a finite tensor category equipped with a (quasi)fiber functor, and set 
$H=\End F$. Then $H$ carries a coproduct $\Delta$ defined as follows: 
$$\Delta:H\to H\otimes H=\End{F\times F};\quad 
T\mapsto J\circ T\circ J^{-1}$$
Moreover, one can define a counit  $\epsilon:H\to\K$ by $\epsilon(T)=T_{\vert F(\textbf1)}$ and an 
antipode $S:H\to H$  by $S(T)_{\vert F(X)}=(T_{\vert F(X^*)})^*$ (in the quasi-case this depends on 
the choice of the identification  $j_X:F(X)^*\to F(X^*)$). \\
This gives $H$ a (quasi-)Hopf algebra structure (the choice of $j_X$ has to do with Drinfeld's 
special elements $\alpha,\beta\in H$). Thus we have bijections: 
{\large
\begin{eqnarray*}
\substack{\textrm{Finite tensor categories with quasifiber} \\
\textrm{functor up to equivalence and changing} \\
\textrm{quasitensor structure of the functor}}  & \longleftrightarrow & 
\substack{\textrm{Finite dimensional quasi-Hopf algebras} \\
\textrm{up to isomorphism and twisting}} \\
\substack{\textrm{Finite tensor categories with fiber} \\
\textrm{functor up to equivalence}} & \longleftrightarrow & 
\substack{\textrm{Finite dimensional Hopf algebras} \\
\textrm{up to isomorphism}}
\end{eqnarray*}}

\subsubsection{Braided and symmetric categories}

Let $\cat$ be a monoidal category with a functorial isomophism 
$\sigma:-\otimes-\to-\otimes^{\textrm{op}}-$, where  $X\otimes^{op}Y:=Y\otimes X$. 

For given objects $V_1,\dots,V_n$ in $\cat$, we consider an expression obtained from 
$V_{i_1}\otimes\cdots\otimes V_{i_n}$ by inserting $\textbf1$'s and parentheses, and where 
$(i_1,\dots,i_n)$ is a permutation of $\{1,\dots,n\}$. To any composition $\varphi$ of $\Phi$'s, 
$\lambda$'s, $\mu$'s, $\sigma$'s and their inverses acting on it, we assign an element of the braid 
group ${\rm B}_n$ as follows: assign $1$ to $\Phi$, $\lambda$ and $\mu$, and the generator $\sigma_k$ of 
${\rm B}_n$ to $\sigma_{V_kV_{k+1}}$. 
\begin{dfn}\emph{
A \emph{braided monoidal category} is a monoidal category as above such that the $\varphi$'s depend 
only on their images in the braid group. 
}\end{dfn}
Again, we have a coherence theorem for braided categories: 
\begin{thm}[\cite{JS}]
The data $(\cat,\otimes,\textbf1,\Phi,\lambda,\mu,\sigma)$ defines a braided category if and only if 
$(\Phi,\alpha)$ satisfy the Hexagon axioms: the diagrams \\
$\xymatrix{
(12)3\ar[r]^\Phi\ar[d]^{\sigma\otimes\textrm{id}} & 1(23)\ar[r]^{\sigma_{1,23}} & (23)1\ar[d]^\Phi \\
(21)3\ar[r]^\Phi & 2(13)\ar[r]^{\textrm{id}\otimes\sigma} & 2(31)}$
$\xymatrix{
(12)3\ar[r]^\Phi\ar[d]^{\sigma^{-1}\otimes\textrm{id}} & 1(23)\ar[r]^{\sigma_{1,23}^{-1}} 
& (23)1\ar[d]^\Phi \\
(21)3\ar[r]^\Phi & 2(13)\ar[r]^{\textrm{id}\otimes\sigma^{-1}} & 2(31)}$ \\
are commutative. 
\end{thm}
\noindent\textbf{Remark. }``2(31)'' is short notation for the 3-functor 
$(V_1,V_2,V_3)\mapsto V_2\otimes(V_3\otimes V_1)$. 

To get the definition of a \emph{symmetric monoidal category}, the reader just has to replace the 
braid group ${\rm B}_n$ by the symmetric group $S_n$ in the definition. To say it in another way, a 
symmetric monoidal category is a braided one for which $\sigma$ satisfies 
$\sigma_{VW}\circ\sigma_{WV}=\textrm{id}_{V\otimes W}$. 
\begin{ex}\emph{
Let $H$ be a quasitriangular bialgebra (resp. Hopf algebra), i.e., a bialgebra (resp. Hopf algebra) 
with an invertible element $R\in H\otimes H$ satisfying $\Delta^{\textrm{op}}(x)=R\Delta(x)R^{-1}$, 
$(\textrm{id}\otimes\Delta)(R)=R^{13}R^{12}$ and $(\Delta\otimes\textrm{id})(R)=R^{13}R^{23}$. 
Then $\rep H$ is a braided  monoidal (resp. rigid monoidal, i.e., tensor) category with braiding 
$\sigma_{VW}:a\otimes b\mapsto R^{21}(b\otimes a)$. Moreover, axioms for $R$ are equivalent to the 
requirement that $\rep H$ is braided (it is not difficult to show that the first equation satisfied 
by $R$ is equivalent to the functoriality of $\sigma$, and the two others are equivalent to the 
Hexagon relations). \\
If $R$ is triangular, i.e., $RR^{21}=1\otimes1$ (in particular if $H$ is cocommutative), then 
$\rep H$ becomes a symmetric monoidal (resp. tensor) category. 
}\end{ex}

\subsubsection{The Drinfeld center}

Tannakian formalism tells us that there is a strong link between finite tensor categories and Hopf 
algebras. So it is natural to ask for a categorification of the notion of the Drinfeld double for 
Hopf algebras. 
\begin{dfn}
\emph{The \emph{Drinfeld center $Z(\cat)$} of a tensor category $\cat$ is a new tensor category whose 
objects are pairs $(X,\Phi)$, where $X\in\obj\cat$ and $\Phi:X\otimes-\to-\otimes X$ is a functorial 
isomorphism such that $\Phi_{Y\otimes Z}=(\textrm{id}\otimes\Phi_Z)\circ(\Phi_Y\otimes\textrm{id})$, 
and with morphisms defined by $\Mor{(X,\Phi)}{(Y,\Psi)}:=\{f\in\Mor XY|\forall Z,
(f\otimes\textrm{id})\circ\Phi_Z=\Psi_Z\circ(\textrm{id}\otimes f)\}$. }
\end{dfn}

\begin{prop}
$Z(\cat)$ is a braided tensor category, which is finite if $\cat$ is. 
\end{prop}
\begin{proof}
See \cite{K} for the proof (the finiteness statement can be found for example in \cite{EO1}). Let us 
just note that the tensor product of objects is given by 
$(X,\Phi)\otimes(Y,\Psi)=(X\otimes Y,\Lambda)$, where 
$\Lambda(Z)=(\Phi(Z)\otimes\textrm{id}_Y)\circ(\textrm{id}_X\otimes\Psi(Z))$, the neutral object by 
$(\textbf1,\textrm{id})$, and the braiding by $\sigma_{(X,\Phi),(Y,\Psi)}=\Phi_Y$.
\end{proof}

\begin{thm}
If $\cat$ is a fusion category over $\C$, then $Z(\cat)$ is also fusion. 
\end{thm}
\begin{proof}
This will be a consequence of a more general statement given in subsection \ref{dualcat}. 
\end{proof}

\noindent\textbf{Remark. }In positive characteristic, $Z(\cat)$ is, in general, not fusion. For 
example, if $\cat=\cat(G,1)$ over $\K=\overline{\mathbb F_p}$, then $Z(\cat)=\rep(\K[G]\ltimes\Fun G)$ 
which is not semisimple if $|G|$ is divisible by $p$. 

~\\

\subsection{Fusion rings}

\subsubsection{Realizability of fusion rings}

Broadly speaking, fusion rings are rings which have the basic properties of Grothendieck rings of 
fusion categories. So let us consider a tensor category $\cat$. 
\begin{enumerate}
\item First, we have seen that if $\cat$ is a tensor category, then $A=\Gr\cat$ is a ring 
which is a free ${\mathbb Z}$-module with a 
distinguished basis $\{X_i\}_{i\in I}$ such that $X_0=\textbf1$ and multiplication (=fusion) rule 
$X_i\cdot X_j=\sum_kN_{ij}^kX_k$, $N_{ij}^k\geq0$ (property 1). \\
\item Second, from the semisimplicity condition we have 
\begin{prop}[see e.g. \cite{ENO}]
If $\cat$ is a semisimple tensor category, then for every simple object $V$ one has $V^*\cong~\!^*V$ 
(so $V\cong V^{**}$). 
\end{prop}
\begin{proof}
The coevaluation map provides an embedding $\textbf1\hookrightarrow V\otimes V^*$. Since the category 
is semisimple, it implies that $V\otimes V^*\cong\textbf1\oplus W$, then there exists a projection 
$p:V\otimes V^*\twoheadrightarrow\textbf1$. But in a rigid category, the only simple object $Y$ such 
that $V\otimes Y$ projects on $\textbf1$ is $^*V$. 
\end{proof}

Thus there exists an involution $*:i\mapsto i^*$ of $I$, defining an antiautomorphism of $A=\Gr\cat$, 
and such that $N_{ij}^0=\delta_{ij^*}$ (property 2). 
\end{enumerate}
\begin{dfn}
\emph{A finite dimensional ring with a basis satisfying  properties 1 and 2 is called a based ring, 
or a \emph{fusion ring}. }
\end{dfn}
One of the basic questions of the theory of fusion categories is 
\begin{pb}
Given a fusion ring $A$, can it be realized as the Grothendieck ring of a fusion category? 
If yes, in how many ways?
\end{pb}
This problem is quite nontrivial, so let us start with a series of examples to illustrate it. 

\subsubsection{Some important examples}

In this subsection we work over $\Bbb C$ unless stated otherwise. 

\begin{ex}\emph{
Consider $A=\Z[G]$ for a finite group $G$, with involution $*:g\mapsto g^{-1}$ being the inversion, 
and the fusion rule being the group law. 
\begin{prop}
The set of realizations of $A$ is $H^{3}(G,{{\bf k}^\times})/Out(G)$. 
\end{prop}
\begin{proof}
Indeed, it is easy to see that the only realizations of $A$ are $\cat(G,\omega)$, and two realizations 
corresponding to 3-cocycles $\omega,\omega'$ are equivalent iff the cohomology classes of 
$\omega,\omega'$ are linked by an automorphism of $G$. Since $G$ acts trivially on its cohomology, we 
get the result. 
\end{proof}
}\end{ex}
\begin{ex}\emph{
Consider fusion ring structures on $A=\Z^2$ (as a $\Z$-module). All such rings are of the form \\
$A_n=<\textbf1,X>$ with $X^2=\textbf1+nX$ and $X^*=X$. 
\begin{thm}[\cite{O2}]\label{2obj}
(i) $A_0$ has two realizations: $\cat(\Z_2,1)$ and $\cat(\Z_2,\omega)$, where $\omega$ is the 
nontrivial element in $H^{3}({\Z_2},{{\bf k}^\times})=\Z_2$. \\
(ii) $A_1$ has two realizations: the fusion category of even highest weight $\aff$-modules at level 3, 
and its Galois image. \\
(iii) For all $n>1$, $A_n$ has no realization. 
\end{thm}
\noindent\textbf{Remark. }The categories in part (ii) of theorem \ref{2obj} are called the 
\emph{Yang-Lee categories} and can also be obtained as quotients of the categories of tilting modules 
over the quantum group $U_q(\sld)$, respectively with $q=e^{\pm i\pi/10}$  and $q=e^{\pm3i\pi/10}$. 
}\end{ex}
\begin{ex}\emph{
Let $B_n$ be the ring generated by $X_0,\dots,X_{n-1}$ and $Y$, satisfying the following relations: 
$Y^2=(n-1)Y+\sum_{i=0}^{n-1}X_i$, $XY=YX=Y$, $Y^*=Y$, $X_iX_j=X_{i+j}$ and $X^*=X_{-1}$ (indices are 
taken $\textrm{mod}~n$). 
\begin{thm}[{\cite[Corollary 7.4]{EGO}}]
$B_n$ is realizable if and only if $q:=n+1$ is a prime power. More precisely, it has three 
realizations for q=3, two when q=4 or 8, and only one for other prime powers. One of the realizations 
is always ${\rm Rep}(\Z_q^\times\ltimes\Z_q)$, the others being obtained by 3-cocycle deformation. 
\end{thm} }\end{ex}
\begin{ex}[Tambara-Yamagami categories, \cite{TY}]\label{ty}\emph{
Let $(G,*)$ be a finite group. Consider $R_G\cong\Z [G]\oplus\Z X$, with fusion product defined by the 
following relations: $X^2=\sum_{g\in G}g$, $gX=Xg=X$, $gh=g*h$, $g^*=g^{-1}$ and $X^*=X$. 
\begin{thm}[\cite{TY}] $R_G$ is realizable if and only if $G$ is abelian. Realizations are 
parametrized by a choice of a sign $\pm$ and a symmetric isomorphism $G\to G^*$ (such an isomorphism 
always exists for abelian groups since it exists for cyclic ones). 
\end{thm}
\noindent If $G=\Z_2$, we obtain the fusion ring corresponding to the \emph{Ising model}: \linebreak 
$R=<\textbf1,g,X>$ with fusion rules $g^2=\textbf1$, $gX=Xg=X$ and $X^2=\textbf1+g$. In this case $R$ 
corresponds to the Grothendieck ring of the category of integrable modules of $\aff$ at level 2 
($V_0=\textbf1$, $V_1=X$ and $V_2=g$). 
}\end{ex}

\subsubsection{The rigidity conjecture}

\begin{cjt}
(i) Any fusion ring  has at most finitely many realizations over ${\bf k}$, up to equivalence (possibly 
none). \\
(ii) The number  of tensor functors between two fixed fusion categories, up to a natural tensor 
isomorphism, is finite. 
\end{cjt}

Thus, the conjecture suggests that fusion categories and functors between them are discrete 
(``rigid'') objects and can't be deformed. It was first proved in the case of unitary categories by 
Ocneanu; thus we call it ``Ocneanu rigidity''. The conjecture is open in general but holds for 
categories over $\C$ (and hence for all fields of characteristic zero). Proving this will be the main 
goal of the next section.  

\section{Ocneanu rigidity}

\subsection{Main results}

\subsubsection{Müger's squared norms}

Let $\cat$ be a fusion category. For every simple object $V$, we are going to define a number 
$|V|^2\in\K^\times$, the \emph{squared norm} of $V$. We have already seen that $V\cong V^{**}$, so 
let us fix an isomorphism $g_V:V\to V^{**}$ and consider its \emph{quantum trace} 
$tr(g_V):=e_{V^*}\circ(g_V\otimes\textrm{id})\circ i_V\in\End{\textbf1}\cong\K$. 

Clearly, this is not an invariant of $V$, since $g_V$ is well defined only up to scaling. However, 
the product $tr(g_V)tr(g_V^{*-1})$ is already independent on the choice of $g_V$ and is an invariant 
of $V$. 

\begin{dfn}[Müger, \cite{Mu1}]
\emph{$|V|^2=tr(g_V)tr(g_V^{*-1})$, and the \emph{global
dimension} of $\cat$ is \footnote{To avoid confusion, we will use the notation
$\textrm{dim}$ for global dimensions, and italic $dim$ for
vector space dimensions.} 
$$\textrm{dim}\cat=\sum_{V~\textrm{simple}}|V|^2\,.$$
If $\textrm{dim}\cat\neq0$, we say that $\cat$ is \emph{nondegenerate}. }
\end{dfn}
\begin{dfn}
\emph{A \emph{pivotal structure} on $\cat$ is an isomorphism of tensor functors $g: Id\to **$. A 
category equipped with a pivotal structure is said to be a {\em pivotal category}. }
\end{dfn}
In a pivotal tensor category, we can define dimensions of objects by $\textrm{dim}V=tr(g_V)$. The 
following obvious properties hold: $\textrm{dim}(V\otimes W)=\textrm{dim}V\textrm{dim}W$ and 
$|V|^2=\textrm{dim}V\textrm{dim}V^*$. 
\begin{dfn}
\emph{We say that a pivotal structure $g$ is {\it spherical} if $\textrm{dim}V=\textrm{dim}V^*$ for 
all simple objects $V$. }
\end{dfn}
\noindent\textbf{Remarks. } 1. It is not known if every fusion category admits a pivotal or spherical 
structure. \\
\indent2. For a simple object $V$ one has $tr(g_V)\neq0$. Indeed, otherwise 
$\textbf1\hookrightarrow V\otimes V^*\twoheadrightarrow\textbf1$, and then the multiplicity 
$[V\otimes V^*:\textbf1]\geq2$, which is impossible in a semisimple category. \\
\begin{ex}\label{kap}
\emph{Let $H$ be a finite dimensional semisimple Hopf algebra over $\K$. Since $\K$ is algebraically 
closed, it is equivalent to saying that $H$ has a decomposition: 
$$H=\bigoplus_{V~\textrm{simple}}\End V\,.$$
It is well-known that the squared antipode $S^2$ is an inner automorphim 
($\exists g\in H^{\times}, S^2(x)=gxg^{-1}$); this is nothing but the statement (proved above) that 
$V$ is isomorphic to $V^{**}$ for simple $H$-modules $V$. Thus $|V|^2=tr_V(S^2_{|\End V})$ and 
$\textrm{dim}(\rep H)=tr_H(S^2)$. \\
It is conjectured (by Kaplansky, \cite{K}) that $S^2=1$; this would imply that $\rep H$ admits a 
spherical structure, such that $|V|^2=dim(V)^2$ and $\textrm{dim}(\rep H)=dim(H)$. For $\K=\C$, this 
is the well-known Larson-Radford theorem \cite{LR}. }
\end{ex}

\subsubsection{Main theorems}

\begin{thm}[Ocneanu, Blanchard-Wassermann, see \cite{BW,ENO}]\label{obw}
If $\cat$ is nondegenerate, then 1) it has no nontrivial first order deformations of its associativity 
constraints, and 2) any tensor functor from $\cat$ has no nontrivial first order deformations of its 
tensor structure. 
\end{thm}
\begin{thm}[\cite{ENO}]
Any fusion category over $\C$ is nondegenerate. 
\end{thm}
The first theorem implies Ocneanu rigidity for nondegenerate fusion categories (see \cite[7.3]{ENO} 
for the precise argument), and the second one proves the rigidity conjecture for fusion categories 
over $\C$. \\
In order to prove these theorems, we have to introduce and discuss the notions of module categories 
and weak Hopf  algebras. 

\subsection{Module categories}

We have seen that the notion of a tensor category is the categorification of the notion of a ring. 
Similarly, the notion of a module category which we are about to define is the categorification of the 
notion of a module over a ring. 

Let $\cat$ be a tensor category. 
\begin{dfn}
\emph{A \emph{left module category over} $\cat$ is an abelian category $\modu$ with an exact 
bifunctor $\otimes:\cat\times\modu\to\modu$ and functorial isomorphisms 
$\alpha:(-\otimes-)\otimes\bullet\to-\otimes(-\otimes\bullet)$ and 
$\eta:\textbf1\otimes\bullet\to\bullet$ (where $\bullet\in \modu$) such that for any two functors 
obtained from $-\otimes\cdots-\otimes\bullet$ by inserting $\textbf1$'s and parenthesis, all 
functorial isomorphisms between them composed of $\Phi^{\pm1}$'s, $\mu^{\pm1}$'s, $\alpha^{\pm1}$'s 
and $\eta^{\pm1}$'s are equal. }
\end{dfn}
The definition of a \emph{right module category over} $\cat$ is analogous. We also leave it to the 
reader to define equivalence of module categories. 

There is an analog of the MacLane coherence theorem for module categories which claims that it is 
sufficient for $\Phi$, $\mu$, $\alpha$ and $\eta$ to make the following diagrams commute: \\
$\xymatrix{
((-\otimes-)\otimes-)\otimes\bullet\ar[r]^{\Phi^{1,2,3}\otimes\textrm{id}}\ar[d]^{\alpha^{12,3,4}} & 
(-\otimes(-\otimes-))\otimes\bullet\ar[r]^{\alpha^{1,23,4}} & 
-\otimes((-\otimes-)\otimes\bullet)\ar[d]^{\textrm{id}\otimes\alpha^{2,3,4}} \\
(-\otimes-)\otimes(-\otimes\bullet)\ar[rr]^{\alpha^{1,2,34}} & & 
-\otimes(-\otimes(-\otimes\bullet))}$ \\ and 
$\xymatrix{
(-\otimes\textbf1)\otimes\bullet\ar[r]^{\alpha_{-,\textbf1,\bullet}}\ar[dr]_{\mu\otimes\textrm{id}} & 
-\otimes(\textbf1\otimes\bullet)\ar[d]^{\textrm{id}\otimes\eta} \\
& -\otimes\bullet}$ \\

\begin{exs}\emph{
(i) $\cat$ is a left module category over itself. \\
\indent(ii) Define the tensor category $\cat^{\textrm{op}}$, which coincides with $\cat$ as an 
abstract category, and has reversed tensor product $\otimes^{\textrm{op}}$, which is defined by 
$X\otimes^{\textrm{op}}Y=Y\otimes X$. The associativity and unit morphisms are defined in an obvious 
manner. Then $\cat$ is a right module category over $\cat^{\textrm{op}}$.  \\
\indent(iii) We deduce from (i) and (ii) that $\cat$ is a left module category over 
$\cat\boxtimes\cat^{\textrm{op}}$. \\
\indent(iv) If $\cat=\vect$ and $\modu=\rep A$ for a given algebra $A$ over $\K$, then $\modu$ is a 
left module category over $\cat$. 
}\end{exs}
Note that if $\modu$ is a left (right) module category over $\cat$, then its Grothendieck group 
$\Gr\modu$ is a left (respectively, right) $\Gr\cat$-module, with a distinguished basis $M_j$ and 
positive structure constants $N_{ij}^r$ such that $X_i\cdot M_j=\sum_rN_{ij}^rM_r$. In this way, we 
can associate  to any object $X\in\obj\cat$ its left (right) multiplication matrix $N_X$, which has 
positive entries, and in the semisimple case $N_{X^*}=N_X^T$. \\

If $\cat$ is a fusion category, we will be interested in semisimple finite module categories over 
$\cat$. Such a module category is called \emph{indecomposable} if $\modu$ is not module equivalent to 
$\modu_1\oplus\modu_2$ for nonzero module categories $\modu_i$, $i=1,2$. 

As was mentioned above, the theory of module categories should be viewed as a categorical analog of 
the theory of modules (representation theory). Thus the main problem in the theory of module 
categories is 
\begin{pb}
Given a fusion category $\cat$, classify all indecomposable module categories over $\cat$ which are 
finite and semisimple. 
\end{pb}
The answer is known only for a few particular cases. For example, one has the following result 
(see \cite{KO,O1} for proof and references): 
\begin{thm}
If $\cat$ is the category of integrable modules over $\aff$ at level $l$, then semisimple finite 
indecomposable module categories over $\cat$ are in one-to-one correspondence with  simply laced 
Dynkin diagrams of ADE type and with Coxeter number $h=l+2$. 
\end{thm}

\subsubsection{The category of bimodules}

Let $\cat$ be a tensor category. A structure of a left module category over $\cat$ on an abelian 
category $\modu$ is the same thing as a tensor functor $\cat\to\textrm{Fun}(\modu,\modu)$ 
($\textrm{Fun}(\modu,\modu)$ is the monoidal category whose objects are exact functors from $\modu$ 
to itself, morphisms are natural transformations, and the tensor product is just the composition of 
functors). This is just the categorical analog of the tautological statement that a module $M$ over a 
ring $A$ is the same thing as a representation $\rho: A\to {\rm End}(M)$. 

If $\modu$ is semisimple and finite, then $\modu\cong\rep A$ as an abelian category for a (nonunique) 
finite dimensional semisimple algebra $A$. Therefore, structures of a left module category over 
$\cat$ on $\modu$ are in one-to-one correspondence with tensor functors 
$\cat\to\textrm{Fun}(\modu,\modu)=A\textrm{-bimod}$. \\
\textbf{Remark. }
In particular, if $\modu$ has only one simple object (i.e. $\modu\cong\vect$ as an abelian category), 
then $\cat$-module category structures on $\modu$ correspond to fiber functors on $\cat$. \\

Let us consider more closely the structure of the category $A$-bimod. Its tensor product 
$\tilde\otimes$ is the tensor product over $A$. The simple objects in this category are 
$M_{ij}=\textrm{Hom}_{\K}(M_i,M_j)$, where $M_i\in\obj\modu$ are simple $A$-modules; and we have 
$M_{ij}\tilde\otimes M_{i'j'}=\delta_{i'j}M_{ij'}$. Thus $A\textrm{-bimod}$ is finite semisimple and 
satisfies all the axioms of a tensor category except one: $\textbf1=\oplus_iM_{ii}$ is not simple, 
but semisimple. 
\begin{dfn}\emph{
A \emph{multitensor category} is a category which satisfies all axioms of a tensor category except 
that the neutral object is only semisimple. \\
A \emph{multifusion category} is a finite semisimple multitensor category. 
}\end{dfn}
Thus, $A\textrm{-bimod}$ is a multifusion category. 

\subsubsection{Construction of module categories over fusion categories}

Let $B$ be an algebra in a fusion category $\cat$. The category $\modu$ of right $B$-modules in 
$\cat$ is a left module category over $\cat$: let $X\in\obj\cat$ and $M$ be a right $B$-module 
($M\otimes B\to M$), then the composition 
$(X\otimes M)\otimes B\tilde\to X\otimes(M\otimes B)\to X\otimes M$ gives us the structure of a right 
$B$-module on $X\otimes M$ (and so it defines a structure of left $\cat$-module category on $\modu$). 
We will consider the situation when $\modu$ is semisimple; in this case the algebra $B$ is said to be 
semisimple. 
\begin{thm}[\cite{O1}]\label{modca}
Any semisimple finite indecomposable module category over a fusion category can be constructed in 
this way (but nonuniquely). 
\end{thm}
\begin{ex}\label{grtheo}\emph{
Let us consider the category $\cat(G,\omega)$, with $G$ a finite group and 
$\omega\in Z^3(G,\K^{\times})$ a 3-cocycle. Let $H\subset G$ be a subgroup such that 
$\omega_{\vert H}=\textrm d\psi$ for a cochain $\psi\in C^2(H,\K^{\times})$. Define the twisted group 
algebra $B=\K_\psi[H]$: $B=\oplus_{h\in H}V_h$ as an object of $\cat$ (where $V_h$ is the 
1-dimensional module corresponding to $h\in H$), and the multiplication map $B\otimes B\to B$ is 
given by $\psi(g,h){\rm Id}: V_g\otimes V_h\to  V_{gh}=V_g\otimes V_h$. The condition 
$\omega_{\vert H}=\textrm d\psi$, which can be rewritten as 
$\psi(h,k)\psi(g,hk)\omega(g,h,k)=\psi(gh,k)\psi(g,h)$ for all $g,h,k\in H$, assures the 
associativity of the product for $B$ (i.e., $B$ is an algebra in $\cat(G,\omega)$). We call 
$\modu(H,\psi)$ the category of right $B$-modules in $\cat(G,\omega)$. 
}\end{ex}
\begin{thm}[\cite{O3}]
Assume ${\rm char}({\bf k})$ does not divide $\vert G\vert$. All semisimple finite indecomposable 
module categories over $\cat(G,\omega)$ have this form. Moreover, two module categories 
$\modu(H_1,\psi_1)$ and $\modu(H_2,\psi_2)$ are equivalent if and only if the pairs $(H_1,\psi_1)$ 
and $(H_2,\psi_2)$ are conjugate under the adjoint action of $G$. 
\end{thm}
\begin{proof}
Let $\modu$ be an indecomposable module category over $\cat(G,\omega)$. Since for every simple object 
we have $X=V_g$, $X\otimes X^*=V_g\otimes V_{g^{-1}}=\textbf1$, the multiplication matrix $N_X$ of $X$ 
satisfies the equation $N_XN_X^T=\textrm{id}$ and thus $N_X$ is a permutation matrix. So we have a 
group homomorphism $G\to\textrm{Perm}(\textrm{simple}(\modu))$. But $\modu$ is indecomposable, 
therefore $G$ acts transitively on $Y:=\textrm{simple}(\modu)$ and so $Y=G/H$. \\
Thus $\modu$ is the category of right $B$-modules in $\cat(G,\omega)$, where $B=\K_\psi[H]$ for a 
$2$-cochain $\psi\in C^2(H,\K^\times)$. 
The associativity condition for the product in $B$, as we saw above, is equivalent to 
$\psi(h,k)\psi(g,hk)\omega(g,h,k)= \psi(gh,k)\psi(g,h)$ (i.e., $\omega_{\vert H}=\textrm d\psi$). We 
are done. 
\end{proof}

\subsection{Weak Hopf algebras}

Tensor functors $\cat\to A\textrm{-bimod}$ are a generalization of fiber functors (which are obtained 
when $A=\K$). So it makes sense to generalize reconstruction theory for them. This leads to 
\emph{Hopf algebroids}, or, in the  semisimple case, to \emph{weak Hopf algebras}. 

\subsubsection{Definition and properties of weak Hopf algebras}

\begin{dfn}[\cite{BNS}]\emph{
A \emph{weak Hopf algebra} is an associative unital algebra $(H,m,1)$ together with a coproduct 
$\Delta$, a counit  $\epsilon$, and an antipode $S$ such that: \\
1) $(H,\Delta,\epsilon)$ is a coassociative counital coalgebra. \\
2) $\Delta$ is a morphism of associative algebras (not necessary unital). \\
3) $(\Delta\otimes\textrm{id})\circ\Delta(1)=(\Delta(1)\otimes1)\cdot(1\otimes\Delta(1))
=(1\otimes\Delta(1))\cdot(\Delta(1)\otimes1)$ \\
4) $\epsilon(fgh)=\epsilon(fg_1)\epsilon(g_2h)=\epsilon(fg_2)\epsilon(g_1h)$ \\
5) $m\circ(\textrm{id}\otimes S)\circ\Delta(h)=
(\epsilon\otimes\textrm{id})\circ(\Delta(1)\cdot(h\otimes1))$ \\
6) $m\circ(S\otimes\textrm{id})\circ\Delta(h)=
(\textrm{id}\otimes\epsilon)\circ((1\otimes h)\cdot\Delta(1))$ \\
7) $S(h)=S(h_1)h_2S(h_3)$
}\end{dfn}

Here we used Sweedler's notation: $\Delta_k(x)=x_1\otimes x_2\otimes...\otimes x_k$ ($\Delta_k$ is the 
$k$-fold coproduct and summation is implicitly assumed). 

\medskip 
\noindent\textbf{Remarks. }1. The notion of finite dimensional weak Hopf algebra is self-dual, i.e., 
if  $(H,m,1,\Delta,\epsilon,S)$ is a finite dimensional weak Hopf algebra then 
$(H^*,\Delta^*,\epsilon^*,m^*,$ $1^*,S^*)$  is also a finite dimensional weak Hopf algebra. \\
2. Let $H$ be a weak Hopf algebra. $H$ is a Hopf algebra if and only if $\Delta(1)=1\otimes1$ (that 
is equivalent to  the requirement that $\epsilon$ is an associative algebra morphism). 

The linear maps $\epsilon_t:h\mapsto\epsilon(1_1h)1_2$ and $\epsilon_s:h\mapsto1_1\epsilon(h1_2)$ 
defined by $5)$ and $6)$ in the definition are called the \emph{target} and 
\emph{source counital maps} respectively. The images $A_t=\epsilon_t(H)$ and $A_s=\epsilon_s(H)$ are 
the target and source \emph{bases} of  $H$. 
\begin{prop}[{\cite[Section 2]{NV}}]\label{bases}
$A_t$ and $A_s$ are semisimple algebras that commute with each other, and $S_{\vert A_t}:A_t\to A_s$ 
is an algebra antihomomorphism. 
\end{prop}
An especially important and tractable class of weak Hopf algebras is that of regular weak Hopf 
algebras, defined as follows. 
\begin{dfn}\emph{
A weak Hopf algebra $H$ is regular if $S^2=\textrm{id}$ on $A_t$ and $A_s$. 
}\end{dfn}

>From now on all weak Hopf algebras we consider will be assumed regular. 

Let $H$ be a finite dimensional weak Hopf algebra and consider the category $\cat=\rep H$. One can 
define the  tensor product $V\otimes W$ of two representations: $V\otimes W:=\Delta(1)(V\otimes_\K W)$ 
as a vector space, and the  action of any $x\in H$ on $V\otimes W$ is given by $\Delta(x)$. As in the 
case of a Hopf algebra, the associativity morphism is the identity, $\epsilon_t$ gives  $A_t$ the 
structure of an $H$-module which is the neutral object in $\cat$, and the antipode $S$ allows us to 
define duality. This endows $\cat$ with the structure of a finite tensor category 
\cite[Section 4]{NTV}. 

In the case when $H$ is regular, each $H$-module $M$ is also an $A_t\otimes A_s$-module (by 
Proposition \ref{bases}), and hence it is an $A_t$-bimodule (since $A_s=A_t^{\textrm{op}}$). Moreover, 
the forgetful functor $\cat=H\textrm{-mod}\to A_t\textrm{-bimod}$ is tensor.

\subsubsection{Reconstruction theory}

Let $\cat$ be a finite tensor category, $A$ a finite dimensional semisimple algebra and 
$F:\cat\to A\textrm{-bimod}$ a tensor functor. Assume that the sizes of the matrix blocks of $A$ are 
not divisible by ${\rm char}(\K)$ (for example, $A$ is commutative or ${\rm char}(\K)=0$). 

Consider $H=\textrm{End}_{\K}(F)=\End{\overline F}$, where $\overline F$ is the composition of $F$ 
with the forgetful functor ${\rm Forget}$ to vector spaces; it is a unital associative algebra. Since 
any $F(X)$ is an $A$-bimodule, there exists an algebra antihomomorphism $s:A\to H$ and an algebra 
morphism $t:A\to H$ such that $[s(a),t(a')]=0$ for all $a,a'\in A$. Moreover, we can define a kind of 
coproduct $\overline\Delta:\textrm{End}_{\K}(F)\to\textrm{End}_{\K}(F\times F)$ in the same way as for 
tannakian formalism: $\overline\Delta(T)=J\circ T\circ J^{-1}$. 
Thus $\overline\Delta(T)$ can be interpreted as an element 
$$\overline\Delta(T)\in H\otimes_AH=H\otimes H/<t(a)x\otimes y-x\otimes s(a)y>,$$
such that $\overline\Delta(T)(t(a)\otimes 1+1\otimes s(a))=0$ for all $a\in H$. Now, since $A$ is 
semisimple, there is a canonical map 
$$\eta:H\otimes_AH\to H\otimes_{\K}H;\quad 
m\otimes n\mapsto\sum_ime_i\otimes e^in$$
for dual bases $(e_i)$ and $(e^i)$ of $A$ relatively to the pairing $(a,b)=tr_A(L_a L_b)$, where 
$L_a$ is the operator of left multiplication by $a$ (note that because of our assumption on the block 
sizes this pairing is nondegenerate). We can thus define the ``true'' coproduct 
$\Delta=\eta\circ\overline\Delta:H\to H\otimes H$ which turns out to be coassociative. \\
One can also define a counit $\epsilon:H\to\K$ by $\epsilon(T)=tr_A(T_{\vert\overline F(\textbf1)})$ 
and an antipode  $S:H\to H$ by $S(T)_{\vert\overline F(X)}=(T_{\vert\overline F(X^*)})^*$. 
\begin{thm}
The associative unital algebra $H$ equipped with $\Delta$, $\epsilon$ and $S$ as above is a regular 
weak Hopf algebra. Moreover, $\cat\cong{\rm Rep}H$ as a tensor category. 
\end{thm}
Thus, given a tensor category $\cat$ over $\K$ and a finite dimensional semisimple algebra $A$ with 
block sizes not divisible by ${\rm char}(\K)$, we have bijections (modulo appropriate equivalences): \\
$
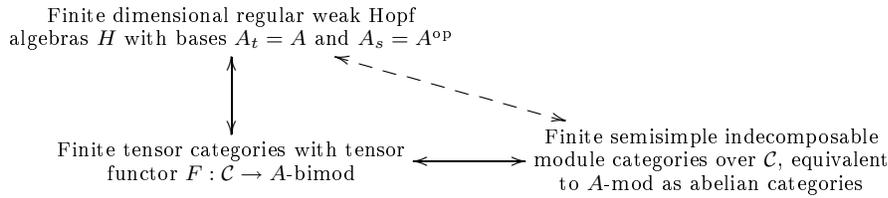
\xymatrix{
\footnotesize{
\txt{Finite dimensional regular weak Hopf\\algebras $H$ with bases $A_t=A$ and $A_s=A^{\textrm{op}}$}}
\ar@{<->}[d]\ar@{<-->}[dr] & \\
\footnotesize{\txt{Finite tensor categories with tensor\\functor $F:\cat\to A\textrm{-bimod}$}}
\ar@{<->}[r] & 
\footnotesize{\txt{Finite semisimple indecomposable\\module categories over $\cat$, equivalent\\
to $A$-mod as abelian  categories}}}$ \\

If $\cat$ is a fusion category, then $\cat$ is a semisimple module category over itself. So 
$\cat\cong\rep H$ as a tensor category for a semisimple weak Hopf algebra $H$ with base 
$A=\oplus_{i\in I}\K_i$. 
\begin{cor}[Hayashi]
Any fusion category is the representation category of a finite dimensional semisimple weak Hopf 
algebra with a commutative base. 
\end{cor}
\noindent{\bf Remark.}
It is not known to us if there exists a (nonsemisimple) finite tensor category which is not the 
category of representations of a weak Hopf algebra (i.e. does not admit a semisimple module category). 
Finding such a category is an interesting open problem. 

\subsection{Proofs}

\subsubsection{Nondegeneracy of fusion categories over $\C$}

\begin{prop}[\cite{N}, \cite{ENO}]
In any fusion category, there exists an isomorphism of tensor functors $\delta:{\rm id}\to ****$. 
\end{prop}
\begin{proof}
Recall that $\cat\cong\rep H$ for a finite dimensional semisimple regular weak Hopf algebra $H$. In 
the semisimple case, the generalization of Radford's $S^4$ formula by Nikshych \cite[Section 5]{N} 
tells us that: 
$$\exists a\in G(H), \forall x\in H, S^4(x)=a^{-1}xa\,,$$
where $a\in G(H)$ means $a$ is invertible and $\Delta(a)=\Delta(1)(a\otimes a)=(a\otimes a)\Delta(1)$ 
(i.e., $a$ is a grouplike element). Thus we can define $\delta$ by $\delta_V=a^{-1}|_V$. Then for 
every $H$-modules $V$ and $W$, the fact that $\delta_{V\otimes W}=\delta_V\otimes\delta_W$ follows 
from the grouplike property of $a$. 
\end{proof}
\begin{thm}[\cite{ENO}]
For fusion categories over $\C$, for any simple object $V$ one has $\vert V\vert^2>0$. 
\end{thm}
In particular, this implies that for any fusion category $\cat$ over $\C$ one has 
$\textrm{dim}\cat\geq1$ and so is nondegenerate. 

{\bf Question.} Does there exist $\epsilon>0$ such that for every fusion category $\cat$ over $\C$ 
which is not $\vect$, $\textrm{dim}\cat>1+\epsilon$?

\begin{proof}[Proof of the theorem]
First do the pivotal case. In this case $\textrm{dim}(V\otimes W)=\textrm{dim}V\textrm{dim}W$ for all 
objects $V,W$, thus $d_id_j=\sum_kN_{ij}^kd_k$, where $d_i=\textrm{dim}(X_i)$ are the dimensions of 
the simple objects. In a shorter way we can rewrite these equalities as $N_i\vec d=d_i\vec d$, where 
$\vec d=(d_0,\dots,d_{n-1})$. \\
For all $i,j,k\in I$, 
\begin{align*}
N_{i^*j}^k & =  
dim(\Mor{X_i^*\otimes X_j}{X_k})
=dim(\Mor{X_j}{X_i\otimes X_k})
\quad\textrm{(by rigidity)} \\
& =  dim(\Mor{X_i\otimes X_k}{X_j})
\quad\textrm{(by semisimplicity)} \\
& =  N_{ik}^j \end{align*}
Therefore $N_i^TN_i\vec d=N_{i^*}N_i\vec d=d_{i^*}d_i\vec d =\vert X_i\vert^2\vec d$, so 
$\vert X_i\vert^2$ is an eigenvalue of $N_i^TN_i$ associated to $\vec d\neq\vec0$ and consequently 
$\vert X_i\vert^2>0$. 

Now we extend the argument to the non-pivotal case. Let us define the \emph{pivotal extension} 
$\overline{\cat}$ of $\cat$, which is the fusion category whose simple objects are pairs $(X,f)$: $X$ 
is simple in $\cat$ and $f:X\tilde\to X^{**}$ satisfies  $f^{**}f=\delta_X$ for the isomorphism of 
tensor functors $\delta:\textrm{id}\to ****$ constructed above. The category $\overline{\cat}$ has a 
canonical pivotal structure $(X,f)\to(X^{**},f^{**})$ (which is given by $f$ itself), thus 
$\vert(X,f)\vert^2>0$. Finally the forgetful functor $\overline{\cat}\to\cat;(X,f)\mapsto X$ preserves 
squared norms, and so $\vert X\vert^2>0$. 
\end{proof}

\subsubsection{Proof of Ocneanu rigidity: the Davydov-Yetter cohomology}

Let $\dat$ be a tensor category. Define the following cochain complex attached to $\dat$: 
\begin{itemize}
\item $C^n(\dat)=\End{T_n}$, where $T_n$ is the $n$-functor 
$\dat^n\to\dat;(X_1,\dots,X_n)\mapsto X_1\otimes...\otimes X_n$ ($T_0=\textbf1$ and 
$T_1=\textrm{id}$). 
\item The differential $\textrm d:C^n(\dat)\to C^{n+1}(\dat)$ is given by 
$$\textrm df=\textrm{id}\otimes f_{2,\dots,n+1}-f_{12,3,\dots,n+1}+\cdots
+(-1)^nf_{1,\dots,n-1,nn+1}+(-1)^{n+1}f_{1,\dots,n}\otimes\textrm{id}$$
\end{itemize}
$H^n(\dat)$ is the $n$-th space of the \emph{Davydov-Yetter cohomology} (\cite{Dav,Y}). 
\begin{ex}\emph{ Assume $\dat=\rep H$ for a Hopf algebra $H$. Then 
$C^n(\dat)=(H^{\otimes n})^{H_{ad}}$. $(C^n,\textrm d)$ is a subcomplex of the co-Hochschild complex 
for $H$ with trivial coefficients. 
}\end{ex}
\begin{prop}[see \cite{Y}]
$H^3(\dat)$ and $H^4(\dat)$ respectively classify first order deformations of associativity 
constraints in $\dat$ and obstructions to these deformations. 
\end{prop}
\begin{exs}\emph{
(i) Let $G$ be a finite group and $\dat=\cat(G,1)$. Then $H^i(\dat)=H^{i}(G,\K)$, and thus $H^i(\dat)=0$ 
for $i>0$ if $\K=\C$ or $\vert G\vert$ and $\textrm{char}(\K)$ are coprime. \\
\indent (ii) Let $G$ be a semisimple complex Lie group with Lie algebra $\al$ and consider 
$\cat=\rep G$. Then $H^i(\cat)=(\wedge^i\al)^G=H^{i}(G,\C)$. In particular, $H^3(\cat)=\C$ and 
$H^4(\cat)=0$. So there exists a unique one-parameter deformation of $\cat=\rep G$ which is realized 
by $\rep U_\hbar(\al)$. 
}\end{exs}
The next result implies in particular the first part of Theorem \ref{obw}. 
\begin{thm}[\cite{ENO}]
Let $\dat$ be a nondegenerate fusion category over $\K$. Then for all $i>0$, $H^i(\dat)=0$. 
\end{thm}
\begin{proof}
The proof is based on the notion of categorical integral. \\
Suppose that $f\in C^n(\dat)$ (for $X_1,\dots,X_n$, 
$f_{X_1,\dots,X_n}:X_1\otimes\cdots\otimes X_n\to X_1\otimes\cdots \otimes X_n$). Define 
$\int f\in C^{n-1}$ in the following way: for $X_1,\dots,X_{n-1}\in\obj\dat$, 
$$(\int f)_{X_1,\dots,X_{n-1}}
=\sum_{V~\textrm{simple}}tr_V((\textrm{id}\otimes g_V)\circ f_{X_1,\dots,X_{n-1},V}) tr(g_V^{*-1})$$
where  $tr_V((\textrm{id}\otimes g_V)\circ f_{X_1,\dots,X_{n-1},V})$ is equal to 
$$(\textrm{id}^{\otimes n}\otimes e_{V^*})\circ
(\textrm{id}^{\otimes(n-1)}\otimes g_V\otimes\textrm{id})
\circ(f_{X_1,\dots,X_{n-1},V}\otimes\textrm{id})
\circ(\textrm{id}^{\otimes n}\otimes i_V)$$
\textbf{Remark. }By definition, $\int\textrm{id}=\textrm{dim}\dat$. \\
Assume now that $f\in Z^{n}(\dat)$ is a cocycle. Then if we put $\varphi=\int f$, we have
\begin{align*}
0 & =  \int{\textrm df} \\
& =  \textrm{id}\otimes\int{f_{2,\dots,n+1}}-\int{f_{12,3,\dots,n+1}}+\cdots \\
& +  (-1)^n\int{f_{1\dots,n-1,nn+1}}+(-1)^{n+1}f_{1,\dots,n}\otimes\int\textrm{id} \\
& = \textrm{id}\otimes\varphi_{2,\dots,n}-\varphi_{12,3,\dots,n}+\cdots \\
& +  (-1)^{n-1}\varphi_{1,\dots,n-1n}+(-1)^n\int{f_{1\dots,n-1,nn+1}}+(-1)^{n+1}
\textrm{dim}\dat\cdot f_{1,\dots,n}
\end{align*}
\begin{lem}[\cite{ENO}]
$\int{f_{1,\dots,n-1,nn+1}}=\varphi_{1,\dots,n-1}\otimes\textrm{id}$. 
\end{lem}
\begin{proof}[Proof of the lemma]
The proof is based the theory of weak Hopf algebras, and we will omit it, see \cite{ENO}, Section 6. 
\end{proof}
Thus when $\textrm{dim}\dat\neq0$, $f=\frac{1}{\textrm{dim}\dat}(-1)^{n-1}\textrm d\varphi$. 
\end{proof}
\noindent\textbf{Remark. }
In the same way, for any tensor functor $F:\cat\to\dat$, one can define a cochain complex 
$C^n_F(\cat)=\End{T_n\circ F^{\otimes n}}$ and a differential 
$\textrm d:C_F^n(\cat)\to C_F^{n+1}(\cat)$ which is given by 
$$\textrm df=\textrm{id}\otimes f_{2,\dots,n+1}-f_{12,3,\dots,n+1}+\cdots
+(-1)^nf_{1,\dots,n-1,nn+1}+(-1)^{n+1}f_{1,\dots,n}\otimes\textrm{id}$$
where $f_{1,\dots,ii+1,\dots,n+1}$ acts on $F(X_1)\otimes\dots\otimes F(X_{n+1})$ as $f$ on 
$F(X_1)\otimes\dots\otimes F(X_i\otimes X_{i+1})\otimes\dots \otimes F(X_{n+1})$ (we have used the 
tensor structure to identify $F(X_i)\otimes F(X_{i+1})$ and $F(X_i\otimes X_{i+1})$). \\
Then one can show (see \cite{ENO}) that the corresponding cohomology spaces $H^i_F(\cat)$ are trivial 
for nondegenerate categories, and that $H^2_F(\cat)$ (resp. $H^3_F(\cat)$) classifies first order 
deformations of the tensor structure of $F$ (resp. obstructions to these deformations). Thus the 
second part of Theorem \ref{obw} is proved. 

\section{Morita theory, modular categories, and lifting theory}

\subsection{Morita theory in the categorical context}

\subsubsection{Dual category with respect to a module category}\label{dualcat}

\begin{pb}
Let $H$ be a finite dimensional (weak) Hopf algebra. $\cat={\rm Rep}(H)$ is a finite tensor category. How to 
describe the category ${\rm Rep}(H^*)$ in terms of $\cat$?
\end{pb}
\noindent The answer is given by the next definitions. 

\begin{dfn}
\emph{A \emph{module functor} between module categories $\modu_1, \modu_2$ over $\cat$ is an additive 
functor $F:\modu_1\to\modu_2$ together with a functorial isomorphism 
$J:F(-\otimes_1\bullet)\to-\otimes_2F(\bullet)$ such that the following diagrams commute: \\
$\xymatrix{
F((-\otimes_\cat-)\otimes_1\bullet)\ar[r]^{F(\alpha)}\ar[d]^J & 
F(-\otimes_1(-\otimes_1\bullet))\ar[r]^J & 
-\otimes_2F(-\otimes_1\bullet)\ar[d]^{\textrm{id}\otimes J} \\
(-\otimes_\cat-)\otimes_2F(\bullet)\ar[rr]^\alpha & & -\otimes_2(-\otimes_2F(\bullet))}$ \\ and 
$\xymatrix{
F(\textbf1\otimes_1\bullet)\ar[r]^{J_{\textbf1,\bullet}}\ar[dr]_{F(\eta_1)} & 
\textbf1\otimes_2F(\bullet)\ar[d]^{\eta_2} \\ 
& F(\bullet)}$ \\ }
\end{dfn}

Let $\cat$ be a tensor category (not necessarily semisimple) and $\modu$ a left module category over 
$\cat$. 
\begin{dfn}\emph{
The \emph{dual category of $\cat$ with respect to $\modu$} is  the category 
$\cat_\modu^*=\textrm{Fun}_\cat(\modu,\modu)$, the category of module  functors from $\modu$ to 
itself with tensor product being the composition of functors. 
}\end{dfn}

Thus the notion of the dual category is the categorification of the notion of the centralizer of an 
algebra in a module. 

Observe that $\cat_\modu^*$ is a monoidal category and $\modu$ is a left module category over it. 
However, $\cat_\modu^*$ is not always rigid. For example, if $\cat=\vect$ and $\modu=A\textrm{-mod}$ 
for a finite dimensional associative algebra $A$ over $\K$, then 
$\cat_\modu^*=\textrm{Fun}_{\vect}(\modu,\modu)=\textrm{Fun}(\modu,\modu)$. This category contains 
the category $A\textrm{-bimod}$ with tensor product $\otimes_A$ which is not exact if $A$ is not 
semisimple (while it must be exact in the rigid case). 

Thus, to insure rigidity of the dual category, we should perhaps restrict ourselves to a subclass of 
module categories. A subclass that turns out to produce a good theory is that of {\it exact module 
categories}. Namely (see \cite{EO1}), a left module category is called \emph{exact} if for any 
projective object $P$ in $\cat$, and any $X\in\obj\modu$, $P\otimes X$ is also projective. Such a 
category is finite if and only if it has finitely many simple objects. In the particular case of a 
fusion category $\cat$, exactness for module categories coincides with semisimplicity. 

\begin{thm}[\cite{EO1}]
If $\cat$ is a finite tensor category and $\modu$ is a finite indecomposable exact left module 
category over it, then $\cat_\modu^*$ is a finite tensor category. 
\end{thm}

\begin{exs}
\emph{(i) If $\cat=\rep H$ and $\modu=\rep A$ for a finite dimensional regular weak Hopf algebra with 
bases $A,A^{\textrm{op}}$, then $\cat_\modu^*=\rep(H^{*\textrm{op}})$. }

\emph{(ii) Let $\cat=\cat(G,\omega)$ and $\modu=\modu(H,\psi)$ be as in example \ref{grtheo}. Then one 
can consider the category of $B$-bimodules $\cat(G,\omega,H,\psi):=\cat_\modu^*$, where 
$B=\K_\psi[H]$ is the twisted group algebra of $H$ in $\cat$. Such categories are called \emph{group 
theoretical}. 
}\end{exs}

Let $\cat$ be a finite tensor category and $\modu$ a finite indecomposable exact left module category 
over $\cat$. Then one can show (\cite{ENO,EO1,O1}) that the following properties hold: 
\begin{enumerate}
\item $(\cat_\modu^*)_\modu^*=\cat$
\item $(\cat\boxtimes\cat_\modu^*)_\modu^*=Z(\cat)$
\item $\cat_\cat^*=\cat^{\textrm{op}}$ (and then $(\cat\boxtimes\cat^{\textrm{op}})_\cat^*=Z(\cat)$ 
by the previous one).
\item If $\modu=B\textrm{-mod}$ for a semisimple algebra $B$ in a fusion category $\cat$, then 
$\cat_\modu^*=B\textrm{-bimod}$. 
\item If $\cat$ is a nondegenerate fusion category, then $\cat_\modu^*$ is also fusion. Moreover, 
$\textrm{dim}\cat_\modu^*=\textrm{dim}\cat$, and thus $\textrm{dim}Z(\cat)=(\textrm{dim}\cat)^2$. 
\end{enumerate}
{\bf Remark.} Note that property (1) is the categorical version of the double centralizer theorem for 
semisimple algebras (saying that the centralizer of the centralizer of $A$ in a module $M$ is $A$ if 
$A$ is a finite dimensional semisimple algebra). Property (2) is the categorical analog of the 
statement that if $A'$ is the centralizer of $A$ in $M$ then the centralizer of $A\otimes A'$ in $M$ 
is the center of $A$. Finally, property (3) is the categorical version of the fact that the 
centralizer of $A$ in $A$ is $A^{op}$.   

\subsubsection{Morita equivalence of finite tensor categories}

By now, all module categories are supposed to be finite and exact. 
\begin{dfn}\emph{
Two finite tensor categories $\cat$ and $\dat$ are \emph{Morita equivalent} if there exists an 
indecomposable (left) module category $\modu$ over $\cat$ such that $\cat_\modu^*=\dat^{\textrm{op}}$. 
In this case we write $\cat\sim_\modu\dat$. 
}\end{dfn}

Obviously, this notion is the categorical analog of Morita equivalence of associative algebras. 
\begin{prop}[Müger, \cite{Mu1,Mu2}]
Morita equivalence of finite tensor categories is an equivalence relation. 
\end{prop}
\begin{proof}
This relation is reflexive since $\cat_\cat^*=\cat^{\textrm{op}}$. \\
To prove the symmetry, assume that $\cat_\modu^*=\dat^{\textrm{op}}$, and define 
$\modu^\vee:=\textrm{Fun}(\modu,\vect)$. This is a left (indecomposable) module category over $\dat$ 
and $\dat_{\modu^\vee}^*=\cat^{\textrm{op}}$. 
Now let us prove transitivity. Suppose $\cat\sim_\modu\dat$ and $\dat\sim_{\mathcal N}\mathcal E$. 
Take $\mathcal P=\textrm{Fun}_\dat(\modu^\vee,\mathcal N)$ (By analogy with ring theory, we could 
denote this category by $\modu\otimes_\dat\mathcal N$.) Then 
$\cat_{\mathcal P}^*=\mathcal E^{\textrm{op}}$. Thus the transitivity condition is verified. 
\end{proof}
\begin{thm}[\cite{Mu1,Mu2}; see also \cite{O3}]\label{Mug}
Let $\cat\sim_\modu\dat$ be a Morita equivalence of finite tensor categories. Then there is a 
bijection between indecomposable left module categories over $\cat$ and $\dat$. It maps $\mathcal N$ 
over $\cat$ to ${\rm Fun}_\cat(\modu,\mathcal N)$ over $\dat$. 
\end{thm}
This, obviously, is the categorical version of the well known characterization of Morita equivalent 
algebras: their categories of modules are equivalent. 
\begin{cor}[\cite{O3}]\label{O3}
Indecomposable left module categories over $\cat(G,\omega,H,\psi)$ are 
$$\modu(H,\psi,H',\psi'):={\rm Fun}_{\cat(G,\omega)}(\modu(H,\psi),\modu(H',\psi'))$$
\end{cor}

\subsubsection{Application to representation theory of groups}

Let $G$ be a finite group and consider the category $\dat=\rep G$. In fact, $\dat=\cat(G,1)_\modu^*$ 
with $\modu=\modu(G,1)=\vect$. Hence, indecomposable $\dat$-module categories are of the form  
$\modu(G,1,H,\psi)=\rep\C_\psi[H]$. \\

Now recall that fiber functors are classified by module categories with only one simple object. In 
our case it corresponds to the case when $\C_\psi[H]$ is simple, which is equivalent to the 
requirement that $\psi$ is a nondegenerate 2-cocycle, in the sense of the following definition. 
\begin{dfn}\emph{
A 2-cocycle $\psi$ on $H$ is \emph{nondegenerate} if $H$ admits a unique projective irreducible 
representation with cocycle $\psi$ of dimension $\sqrt{\vert H\vert}$. \\
A group $H$ which admits a nondegenerate cocycle is said to be \emph{of central type}. 
}\end{dfn}
\noindent\textbf{Remarks. } 1. It is obvious that a group of central type has order $N^2$, where $N$ 
is an integer. \\
\indent2. Howlett and Isaacs \cite{HI} showed that any group of central type is solvable. This is a 
deep result based on the classification of finite simple groups. 

\begin{thm}[\cite{EG,M}]
Fiber functors on ${\rm Rep}G$ (i.e., Hopf twists on $\C[G]$ up to a gauge) are in one-to-one 
correspondence with pairs $(H,\psi)$, where $H$ is a subgroup of $G$ and $\psi$ a nondegenerate 
2-cocycle on $H$ modulo coboundaries and inner automorphisms. 
\end{thm}
\begin{proof} The theorem follows from Theorem \ref{Mug} and Corollary \ref{O3}. We leave the proof 
to the reader. 
\end{proof}
\begin{cor}[\cite{TY}]
Let $D_8$ be the group of symmetries of the square and $Q_8$ the quaternion group. Then 
${\rm Rep}D_8$ and ${\rm Rep}Q_8$ are not equivalent (although they have the same Grothendieck ring). 
\end{cor}
\begin{proof}[Proof of the corollary]
In $Q_8$, all subgroups of order 4 are cyclic and hence do not admit any nondegenerate 2-cocycle. \\
On the other hand, $D_8$ has two subgroups isomorphic to $\Z_2\times\Z_2$ (not conjugate) and each has 
one nondegenerate 2-cocycle. Thus $Q_8$ has fewer fiber functors (in fact only 1) than $D_8$ (which 
has 3 such). 
\end{proof}
So, we see that one can sometimes establish that two fusion categories are not equivalent (as tensor 
categories) by counting fiber functors. Similarly, one can sometimes show that two fusion categories 
are not Morita equivalent by counting all indecomposable module categories over them (since we have 
seen that Morita equivalent fusion categories have the same number of indecomposable module 
categories). Let us illustrate it with the following example. 

\begin{ex}\emph{ We want to show that $\rep(\Z_p\times\Z_p)$ and $\rep\Z_{p^2}$ are not Morita 
equivalent. \\
\indent First remember that $\rep G=\cat(G,1,G,1)$ and module categories over it are parametrized by 
$(H,\psi)$, where $H$  is a subgroup of $G$ and $\psi\in H^{2}(H,\C^{\times})$. \\
On the one hand, $\Z_{p^2}$ has three subgroups ($\Z_{p^2}$ itself, $\Z_p$, and $1$), all with a 
trivial second cohomology. Thus $\rep\Z_{p^2}$ has 3 indecomposable module categories. On the other 
hand, $\Z_p\times\Z_p$  has $p+3$ subgroups: $\Z_p\times\Z_p$, $p+1$ copies of $\Z_p$, and $1$. 
Moreover, $\Z_p\times\Z_p$ has $p$ 2-cocycles up to coboundaries. Thus $\rep(\Z_p\times\Z_p)$ has 
$2p+2>3$ module categories. $\Box$ }\end{ex}  

\subsection{Modular categories and the Verlinde formula}

Let $\cat$ be a braided tensor category. 
Then we have a canonical (non-tensor) functorial isomorphism $u: \textrm{id}\to **$ given by the 
composition 
$$V\to V\otimes V^*\otimes V^{**}\to V^*\otimes V\otimes V^{**}\to V^{**}$$
(the maps are the coevaluation, the braiding, and the evaluation). This isomorphism is called the 
{\em Drinfeld isomorphism}. Using the Drinfeld isomorphism, we can define a tensor isomorphism 
$\delta: {\rm id}\to ****$  by the formula $\delta_V=(u_{V^*}^*)^{-1}u_V$. 

\begin{dfn}\emph{
A \emph{ribbon category} is a braided tensor category together with a pivotal structure 
$g:\textrm{id}\to**$, such that $g^{**}g=\delta$. 
}\end{dfn}
We refer the reader who wants to learn more about ribbon categories (especially the graphical 
calculus for morphisms, using tangles) to \cite{K}, \cite{BK} or \cite{T}. 

Assume now that $\cat$ is a ribbon category. Recall for any simple object $V\in\cat$ one can define 
the \emph{dimension} $\textrm{dim}V$. It is known (see e.g. \cite{K}) that 
$\textrm{dim}V^*=\textrm{dim}V$. 

For any two objects $V,W$, one can define the number $S_{VW}\in\End{\textbf1}\cong\K$ to be 
{\small
$$(e_{V^*}\otimes e_{W^*})\circ(g_V\otimes\textrm{id}_{V^*}\otimes g_W\otimes\textrm{id}_{W^*})\circ
(\textrm{id}_V\otimes\sigma_{WV^*}\otimes\textrm{id}_{W^*})
\circ(\textrm{id}_V\otimes\sigma_{V^*W}\otimes\textrm{id}_{W^*})\circ(i_V\otimes i_W)\,.$$}
Now assume that $\cat$ is fusion, with simple objects $X_i$'s. 
Then we can define a matrix $S$ with entries $S_{ij}=S_{X_iX_j}$. $S$ has the following properties: 
\begin{enumerate}
\item $S_{ij}=S_{ji}$
\item $S_{ij}=S_{i^*j^*}$
\item $S_{i0}=\textrm{dim}X_i\neq0$
\end{enumerate}
\begin{dfn}\emph{
A ribbon fusion category is called \emph{modular} if $S$ is nondegenerate. 
}\end{dfn}
\begin{prop}[\cite{Mu2,T}]
If $\cat$ is a nondegenerate fusion category with a spherical structure, then $Z(\cat)$ is a modular 
category. 
\end{prop}
\begin{prop}[{\cite[Theorem 3.1.7]{BK}}]
In a modular category $\cat$, 
$$\sum_k  S_{ik}S_{kj}=({\rm dim}\cat)\delta_{ij^*}$$
\end{prop}
Thus if $\cat$ is a modular category, then $\textrm{dim}\cat\neq0$ and we can define new numbers 
$s_{ij}=S_{ij}/\sqrt{\textrm{dim}\cat}$ (here we must make a choice of the square root). 
\begin{thm}[Verlinde formula, \cite{BK}]
$$\sum_\alpha N_{ij}^\alpha s_{\alpha r}=\frac{s_{ir}s_{jr}}{s_{0r}}$$
\end{thm}

So $s_{ir}/s_{0r}$ are eigenvalues of the multiplication matrix $N_i$. In particular, they are 
algebraic integers (i.e. roots of a monic polynomial with integer coefficients - the characteristic 
polynomial of $N_i$). Hence: 
\begin{prop}
For every $r$, $\frac{{\rm dim}\cat}{({\rm dim}X_r)^2}=\frac{s_{ir}s_{i^*r}}{s_{0r}^2}$ is an 
algebraic integer. 
\end{prop}

This result will be very useful to prove classification theorems in section 4. 

\subsubsection{Galois property of the S-matrix}

A remarkable result due to J.~de~Boere, J.~Goeree, A.~Coste and T.~Gannon states that the entries of 
the S-matrix of a modular category lie in a cyclotomic field, see \cite{dBG,CG}. Namely, 
one has the following theorem. 
\begin{thm}\label{app}
Let $S=(s_{ij})_{i,j\in I}$ be the S-matrix of a modular category ${\mathcal C}$. There exists 
a root of unity $\xi$ such that $s_{ij}\in {\mathbb Q}(\xi)$. 
\end{thm}
\begin{proof}
Let $\{ X_i\}_{i\in I}$ be the representatives of isomorphism classes of simple objects of 
${\mathcal C}$; let $0\in I$ be such that $X_0$ is the unit object of ${\mathcal C}$ and the 
involution $i\mapsto i^*$ of $I$ be defined by $X_i^*\cong X_{i^*}$. 
By the definition of modularity, any homomorphism $f: K({\mathcal C}) \to {\mathbb C}$ is of the form 
$f([X_i])=s_{ij}/s_{0j}$ for some well defined $j\in I$. Hence for any automorphism $g$ of 
$\overline{\mathbb Q}$ one has $g(s_{ij}/s_{0j})= s_{ig(j)}/s_{0g(j)}$ for a well defined action of $g$ on $I$. 

Now remember from the previous subsection that one has the following properties: 
$\sum_ks_{ik}s_{kj}=\delta_{ij^*}$, $s_{ij}=s_{ji}$, and $s_{0i^*}=s_{0i}\ne 0$. \\
Thus, $\sum_js_{ij}s_{ji^*}=1$ and hence $(1/s_{0i})^2=\sum_j(s_{ji}/s_{0i}) (s_{ji^*}/s_{0i^*})$. 
Applying the automorphism $g$ to this equation we get
$$g(\frac{1}{s_{0i}^2})=g(\sum_j\frac{s_{ji}}{s_{0i}}\frac{s_{ji^*}}{s_{0i^*}})=
\sum_j\frac{s_{jg(i)}}{s_{0g(i)}}\frac{s_{jg(i^*)}}{s_{0g(i^*)}}=
\frac{\delta_{g(i)g(i^*)^*}}{s_{0g(i)}s_{0g(i^*)}}\,.$$
It follows that $g(i^*)=g(i)^*$ and $g((s_{0i})^2)=(s_{0g(i)})^2$. Hence 
$$g((s_{ij})^2)=g((s_{ij}/s_{0j})^2\cdot (s_{0j})^2)=(s_{ig(j)})^2\,.$$
Thus $g(s_{ij})=\pm s_{ig(j)}$. Moreover the sign $\epsilon_g(i)=\pm 1$ such that 
$g(s_{0i})=\epsilon_g(i)s_{0g(i)}$ is well defined since $s_{0i}\ne 0$, and 
$g(s_{ij})=g((s_{ij}/s_{0j})s_{0j})=\epsilon_g(j)s_{ig(j)}=\epsilon_g(i) s_{g(i)j}$. In particular, 
the extension $L$ of ${\mathbb Q}$ generated by all entries $s_{ij}$ is finite and normal, that is 
Galois extension. 
Now let $h$ be another automorphism of $\overline{\mathbb Q}$. We have 
$$gh(s_{ij})=g(\epsilon_h(j)s_{ih(j)})=\epsilon_g(i)\epsilon_h(j)s_{g(i)h(j)}$$
and 
$$hg(s_{ij})=h(\epsilon_g(i)s_{g(i)j})=\epsilon_h(j)\epsilon_g(i)s_{g(i)h(j)}= gh(s_{ij})$$
and the Galois group of $L$ over ${\mathbb Q}$ is abelian. Now the Kronecker-Weber theorem 
(see e.g. \cite{Ca}) implies the result. 
\end{proof}

\subsection{Lifting theory}

First recall that a fusion category over an algebraically closed field $\K$ can be regarded as a 
collection of finite dimensional vector spaces $H_{ij}^k$ (=${\rm Hom}(X_i\otimes X_j,X_k)$), 
together with linear maps between direct sums of tensor products of these spaces which satisfy some 
equations (given by axioms of tensor categories). Thus one can define a fusion category over any 
commutative ring with $R$ to be a collection of free finite rank $R$-modules $H_{ij}^k$ together with 
module homomorphisms between direct sums of tensor products of them which satisfy the same equations. \\ 
By a realization of a fusion ring $A$ over $R$ we will mean a fusion category over $R$ such that 
$N_{ij}^k:=dim(H_{ij}^k)$ are the structure constants of $A$. 

If $I$ is an ideal in $R$ and $\cat$ a fusion category over $R$ then it is clear how to define the 
reduced (=quotient) fusion category $\cat/I$ over $R/I$ with the same Grothendieck ring. 

Tensor functors between fusion categories over $\K$ can be defined in similar terms, as collections 
of linear maps satisfying algebraic equations; this allows one to define tensor functors between 
fusion categories over $R$ (and their reduction modulo ideals) in an obvious way. 

Now let $\K$ be any algebraically closed field of characteristic $p>0$, $W(\K)$ the ring of Witt 
vectors of $\K$, $I$ the  maximal ideal of $W(\K)$ generated by $p$, and $\mathbb K$ the algebraic 
closure of the fraction field of $W(\K)$ (char$(\mathbb K)=0$). 
\begin{dfn}\emph{
Let $\cat$ be a fusion category over $\K$. A \emph{lifting $\widetilde\cat$ of $\cat$ to $W({\bf k})$} is 
a realization of $\Gr\cat$ over the ring $W(\K)$ together with an equivalence of tensor categories 
$\widetilde\cat/I\tilde\to\cat$. \\
In a similar way, one defines a \emph{lifting of a tensor functor} $F:\cat\to\dat$: it is a tensor 
functor $\widetilde F:\widetilde\cat\to\widetilde\dat$ over $W(\K)$ together with an equivalence of 
tensor functors $\widetilde F/I\tilde\to F$. 
}\end{dfn}
\begin{thm}[\cite{ENO}]
Let $\cat$ be a nondegenerate fusion category over ${\bf k}$. Then there exists a unique lifting of $\cat$ 
to $W({\bf k})$. 
\end{thm}
\begin{proof}
This follows from the fact that liftings are classified by $H^3(\cat)$ and obstructions by 
$H^4(\cat)$. And we know from Section 2 that the Davydov-Yetter cohomology vanishes for nondegenerate 
categories. 
\end{proof}
\begin{thm}[\cite{ENO}]
Let $F:\cat\to\dat$ be a tensor functor between nondegenerate fusion categories over ${\bf k}$. Then there 
exists a unique lifting of $F$ to $W({\bf k})$. 
\end{thm}\label{liffun}
\begin{proof}
Again, liftings of $F$ are parametrized by $H^2_F(\cat)$ and obstructions by $H^3_F(\cat)$, which are 
trivial in the  nondegenerate case. 
\end{proof}
\begin{cor}[\cite{EG2}]
Any semisimple Hopf algebra $H$ over ${\bf k}$ with $tr(S^2)\neq0$ (i.e., also cosemisimple) lifts to 
$\widetilde H$ over  $W({\bf k})$. 
\end{cor}
Hence one can define $\widehat H=\widetilde H\otimes_{W(\K)}\mathbb K$, which is a Hopf algebra over 
a field of charactristic zero. This allows one to extend results from the characteristic zero case to 
positive characteristic. For example, applying the Larson-Radford theorem \cite{LR} (see Corollary 
\ref{lr} below) to $\widehat H$, one can find: 
\begin{cor}[Kaplansky 7th conjecture, \cite{EG2}]
If $H$ is a semisimple and cosemisimple Hopf algebra over any algebraically closed field, then $S^2=1$. 
\end{cor}
\begin{cor}[\cite{ENO}]\label{nbs}
A nondegenerate braided (resp. symmetric) fusion category over ${\bf k}$ is uniquely liftable to a braided 
(resp. symmetric) fusion category over $W({\bf k})$. 
\end{cor}
\begin{proof}A braiding on $\cat$ is the same as a splitting $\cat\to Z(\cat)$ of the natural 
(forgetful) tensor functor $Z(\cat)\to\cat$. Theorem \ref{liffun} implies that such a splitting is 
uniquely liftable. Thus a braiding is uniquely liftable. 

Now prove the result in the symmetric case. A braiding gives rise to a categorical equivalence 
$B:\cat\to\cat^{\textrm{op}}$, and it is symmetric if and only if the composition of $B$ and 
$B^{21}$ is the identity. Hence the corollary follows from Theorem \ref{liffun}. 
\end{proof}

We conclude the section with mentioning a remarkable theorem of Deligne on the classification of 
symmetric fusion categories over $\Bbb C$. 
\begin{thm}[\cite{De}]
Any symmetric fusion category over $\C$ is ${\rm Rep} G$ for a finite group $G$. 
\end{thm}
\noindent With some work, one can extend this result using corollary \ref{nbs}: 
\begin{cor}[\cite{EG3}]
Any symmetric nondegenerate fusion category over ${\bf k}$ (of characteristic $p$) is 
${\rm Rep} G$ for a finite group $G$ of order not divisible by $p$. 
\end{cor}

\section{Frobenius-Perron dimension}

\subsection{Definition and properties}

Let $\cat$ be a finite tensor category with simple objects $X_0,\dots,X_{n-1}$. Then for every object 
$X\in\obj\cat$, we have a matrix $N_X$ of left multiplication by $X$: $[X\otimes X_i:X_j]=(N_X)_{ij}$. 
This matrix has nonnegative entries, and in the Grothendieck ring we have : 
$XX_i=\sum_j(N_X)_{ij}X_j$. 

Let us now recall the classical
\begin{thm}[Frobenius-Perron]
Let $A$ be a square matrix with nonnegative entries. Then 
\begin{enumerate}
\item $A$ has a nonnegative real eigenvalue. The largest such eigenvalue $\lambda(A)$ dominates in 
absolute value all other eigenvalues of $A$. Thus the largest nonnegative eigenvalue of $A$ coincides 
with the spectral radius of $A$. 
\item If $A$ has strictly positive entries, then $\lambda(A)$ is a simple eigenvalue, which is 
strictly positive, and its eigenvector can be normalized to have strictly positive entries. 
Moreover, if $v$ is an eigenvector with strictly positive entries, then the corresponding eigenvalue 
is $\lambda(A)$. 
\end{enumerate}
\end{thm}

Thus to all $X\in\obj\cat$ one can associate a nonnegative number $d_+(X)=\lambda(N_X)$, its 
\emph{Frobenius-Perron dimension}. \\
\begin{exs}\emph{
(i) The Yang-Lee category: $X^2=\textbf1+X$, so $N_X=(\substack{0~1 \\ 1~1})$ and 
$d_+(X)=\frac{1+\sqrt5}{2}$. \\
\indent (ii) Let $\cat=\rep H$ for a finite dimensional quasi-Hopf algebra $H$, then $d_+(X)=dim(X)$ 
for all $H$-modules  $X$. 
}\end{exs}
The following proposition follows from the interpretation of $d_+(X)$ as the spectral radius of $N_X$. 
\begin{prop}
For all objects $X$ of $\cat$, $\frac{\log ({\rm length}(X^{\otimes n}))}{\log n}\to d_+(X)$ when 
$n$ goes to infinity. 
\end{prop}
\begin{thm}[\cite{ENO}, \cite{E}]
The assignment $X\mapsto d_+(X)$ extends to a ring homomorphism 
${\rm Gr}(\cat)\to\RR$. Moreover, $d_+(X_i)>0$ for $i=0,\dots,n-1$. 
\end{thm}
\begin{proof}
Consider $X=\sum_iX_i\in\Gr\cat$ and denote by $M_X$ the matrix of right multiplication by $X$. For 
$i,j\in I$, 
\begin{align*}
(M_X)_{ij} & =  [X_i\otimes X:X_j]\geq dim(\Mor{X_i\otimes X}{X_j}) \\
& = \sum_k dim(\Mor{X_k}{^*\!\!X_i\otimes X_j})>0\,.
\end{align*}
Hence by the Frobenius-Perron theorem, there exists a unique eigenvector of $M_X$ (up to scaling) 
with strictly positive entries, say $R=\sum_i\alpha_iX_i$: $RX=\mu R$ with $\mu=\lambda(M_X)$. Now 
for all $Y\in\Gr\cat$, $(YR)X=\mu YR$ and then by the uniqueness of $R$ there is $\beta_Y\in\RR$ such 
that $YR=\beta_YR$. Since $R$ has positive coefficients, applying again the Frobenius-Perron 
theorem, we obtain $\beta_Y=\lambda(N_Y)=d_+(Y)$. \\
Consequently, $d_+(Y+Z)R=(Y+Z)R=YR+ZR=(d_+(Y)+d_+(Z))R$ and $d_+(YZ)R=YZR=Yd_+(Z)R=d_+(Y)d_+(Z)R$. So 
$Y\mapsto d_+(Y)$ extends to a ring homomorphism $\Gr\cat\to\RR$. 

Suppose $d_+(X_i)=0$, then $X_iR=0$ and hence $X_iX_j=0$ for all $j\in I$, which is not possible. 
Thus $d_+(X_i)>0$. 
\end{proof}

\noindent\textbf{Remark.}
It is clear that the Frobenius-Perron dimension can be defined for any finite dimensional ring with 
distinguished basis and nonnegative structure constants (even if it has no realization) and does not 
depend on the corresponding category. 

\begin{prop}
$d_+$ is the unique character of ${\rm Gr}(\cat)$ that maps elements of the basis to strictly positive 
numbers. 
\end{prop}
\begin{proof}
Let $\chi$ be another such character. Then $\chi(X_i)\chi(X_j)=\sum N_{ij}^k\chi(X_k)$. Thus the 
vector with positive entries $\chi(X_k)$ is an eigenvector of the matrix $N_i$ with eigenvalue 
$\chi(X_i)$. So by the Frobenius-Perron theorem, $\chi(X_i)=d_+(X_i)$. 
\end{proof}

\begin{cor}
Quasitensor functors between finite tensor categories preserve Frobenius\nolinebreak-Perron dimension. 
\end{cor}
\begin{cor}
$d_+(X)=d_+(X^*)$. 
\end{cor}

\subsubsection*{Properties of the Frobenius-Perron dimension}

\begin{enumerate}
\item $\alpha=d_+(X)$ is an algebraic integer (it is a root of the characteristic polynomial of 
$N_X$). 
\item $\forall g\in\textrm{Gal}(\overline\Q/\Q), \vert g\alpha\vert\leq\alpha$ (use part two of the 
Frobenius-Perron theorem). In  particular, $\alpha\geq1$. 
\item $\alpha=1\Leftrightarrow X\otimes X^*=\textbf1$ (in this case $X$ is called \emph{invertible}). 
\begin{proof}
If $X\otimes X^*=\textbf1$, then $1=d_+(\textbf1)=d_+(X)d_+(X^*)$. Since  $d_+(X)\geq1$ and 
$d_+(X^*)\geq1$, we find that $d_+(X^*)=1$. Conversely, consider 
$i_X:\textbf1\hookrightarrow X\otimes X^*$ and compute 
$$d_+(X\otimes X^*)=d_+(\textbf1)+d_+(\textrm{coker}i_X)=1+d_+(\textrm{coker}i_X)\,.$$
Now if $d_+(X)=1$, then $d_+(X\otimes X^*)=1$, so $d_+(\textrm{coker }i_X)=0$ and hence 
$\textrm{coker }i_X\cong0$. Consequently, $i_X$ is an isomorphism and thus 
$\textbf1\cong X\otimes X^*$. 
\end{proof}
\item (\cite{GHJ}) If $\alpha<2$, then $\alpha=2\cos{\frac{\pi}n}$ for $n\geq3$. 
\begin{proof}
Since $d_+$ is a character, $\alpha$ is the largest characteristic value of $N_X$. But the largest 
characteristic value of a positive integer matrix $A$ (i.e., the spectral radius of $\sqrt{AA^T}$) 
is, by Kronecker's theorem, of the form  $2\cos(\frac{\pi}{n})$, or is $\ge 2$. 
\end{proof}
\end{enumerate}
\begin{thm}[\cite{EO1}]\label{qhthm}
Let $\cat$ be a finite tensor category. $\cat\cong{\rm Rep}H$ as a tensor category for a finite 
dimensional quasi-Hopf algebra $H$ if and only if every object $X$ of $\cat$ has an integer 
Frobenius-Perron dimension. 
\end{thm}
\begin{proof}
First suppose that every object $X$ is such that $d_+(X)\in\N$. Then one can consider the object 
$P=\sum_id_+(X_i)P_i$, where $P_i$ are projective covers of $X_i$, and define a functor 
$F:\cat\to\vect;X\mapsto\Mor PX$, which is exact. Since $F(-)\otimes F(-)$ and $F(-\otimes-)$ extend 
to exact functors $\cat\boxtimes\cat\to\vect$ that map simple objects $X_i\boxtimes X_j$ to the same 
images, they are isomorphic. Thus $F$ is quasitensor and $\cat\cong\rep H$. 

If $\cat\cong\rep H$, then reconstruction theory says there exists a quasifiber functor on $\cat$. 
We know that such a functor  preserves Frobenius-Perron dimensions, so they are integers. 
\end{proof}
\begin{cor}
If $H_1$, $H_2$ are finite dimensional quasi-Hopf algebras such that ${\rm Rep}H_1\cong{\rm Rep}H_2$ as 
tensor categories, then  $H_1$ and $H_2$ are equivalent by a twist. 
\end{cor}
\begin{proof}
In the proof of Theorem \ref{qhthm}, there is no choice in the definition of the quasifiber functor 
$F$. Thus (by reconstruction  theory) $H$ is unique up to a twist. 
\end{proof}
\noindent\textbf{Remark. }
This is not true in the infinite dimensional case. For example, consider the category 
$\cat=\rep(SL_q(2))$ of representations of the quantum group $SL_q(2)$ with $q$ not equal to a 
nontrivial root of unity. Then there are many fiber functors on $\cat$ which are not isomorphic 
(even as usual functors). More precisely, for every $m\geq2$ one can find a tensor functor 
$F:\cat\to\vect$ such that $dim(F(V_1))=m$ (where $V_1$ is the standard 2-dimensional representation 
of $SL_q(2)$). Such $F$ can be classified and yield \emph{quantum groups of a non-degenerate bilinear 
form} \cite{B,EO2}. 

Finally, let us give a number-theoretic property of the Frobenius-Perron dimension in a fusion 
category, which  allows one to dismiss many fusion rings as non-realizable. 
\begin{thm}[\cite{ENO}]
If $\cat$ is a fusion category over $\C$, then there exists a root of unity $\xi$ such that for 
every object $X$ of $\cat$  $d_+(X)\in\Z[\xi]$. 
\end{thm}
\begin{ex}\emph{
Consider the fusion ring $A$ with basis $\textbf1,X,Y$ and fusion rules $XY=2X+Y$, $X^2=\textbf1+2Y$ 
and $Y^2=\textbf1+X+2Y$. The computation of $d_+(X)$ reduces to a cubic equation whose Galois group 
is $S_3$. So we cannot find any root of unity $\xi$ such that $d_+(X)\in\Z[\xi]$, and consequently 
$A$ is not realizable. 
}\end{ex}

\subsection{FP-dimension of the category}

Let $\cat$ be a finite tensor category with simple objects $X_0,\dots,X_{n-1}$. We denote by $P_i$ 
the projective cover of $X_i$ ($i=0,\dots,n-1$). 
\begin{dfn}\emph{
The \emph{Frobenius-Perron dimension of the category} $\cat$ is $d_+(\cat)=\sum_id_+(X_i)d_+(P_i)$. 
}\end{dfn}
\begin{exs}\emph{
(i) If $\cat$ is semisimple (and hence fusion), then $d_+(\cat)=\sum_id_+(X_i)^2$. \\
\indent (ii) If $\cat=\rep H$ for a finite dimensional quasi-Hopf algebra $H$, then 
$d_+(\cat)=dim(H)$. 
}\end{exs}

The usefulness of this notion is demonstrated, for example, by the following result. 
\begin{prop}[\cite{EO1}]
The Frobenius-Perron dimension of the category is invariant under Morita equivalence. 
\end{prop}
\noindent Remember that $Z(\cat)$ is Morita equivalent to $\cat\boxtimes\cat^{\textrm{op}}$. Thus we 
have 
\begin{cor}
Let $\cat$ be a finite tensor category. Then $d_+(Z(\cat))=d_+(\cat)^2$. 
\end{cor}
We note that for spherical categories these results appear in \cite{Mu1}, \cite{Mu2}. 

The following theorem plays a crucial role in classification of tensor categories, and in particular 
allows one to show that many fusion rings are non-realizable. 
\begin{thm}[\cite{EO1}]\label{lag}
If $\cat$ is a full tensor subcategory of a finite tensor category $\dat$, then 
$\frac{d_+(\dat)}{d_+(\cat)}$ is an algebraic integer. 
\end{thm}
\begin{exs}\emph{
(i) Let $\dat=\cat(G,1)$ and $\cat=\cat(H,1)$ for a finite group $G$ and its subgroup $H$. Then 
Theorem \ref{lag} says that $\vert H\vert$ divides $\vert G\vert$ (because an algebraic integer 
which is also a rational number is an integer). Thus Theorem \ref{lag} is a categorical generalization 
of Lagrange's theorem for finite groups.\\
\indent (ii) Let $\dat=\rep A$ and $\cat=\rep B$ for a finite dimensional Hopf algebra $A$ and a 
quotient $B=A/I$ of $A$ by a Hopf ideal $I$. Theorem \ref{lag} says $dim(B)$ divides $dim(A)$ (this 
is the famous Nichols-Zoeller theorem \cite{NZ}). The same applies to quasi-Hopf algebras (in which 
case the result is due to Schauenburg, \cite{S}). 
}\end{exs}
\begin{thm}[\cite{ENO}]
If $\cat$ is a fusion category with integer $d_+(\cat)$, then $d_+(X_i)^2\in\N$ for all $i\in I$. 
\end{thm}
\begin{proof}
Let $\cat_{\textrm{ad}}$ be the full tensor subcategory of $\cat$ generated by direct summands of 
$X_i\otimes X_i^*$ ($i\in I$), and define $B=\oplus_i(X_i\otimes X_i^*)$. This object has an integer 
FP dimension: $d_+(B)=d_+(\cat)\in\N$. Then consider $M=N_{B^{\otimes m}}$, the left multiplication 
matrix by $B^{\otimes m}$ in $\cat_{\textrm{ad}}$. This matrix has positive entries for large enough 
$m$ (since any simple object of $\cat_{\textrm{ad}}$ is contained in $B^{\otimes m}$). \\
Let $Y_0,...,Y_p$ be the simple objects of $\cat_{\textrm{ad}}$. The vector 
$(d_+(Y_0),\dots,d_+(Y_p))$ is an eigenvector of $M$ with integer eigenvalue $d_+(B)^n$. 
By the Frobenius-Perron theorem, this eigenvalue is simple. Thus the entries of the eigenvector are 
rational (as $d_+(Y_0)=1$) and hence integer (as they are algebraic integers). Consequently, 
$d_+(X_i\otimes X_i^*)=d_+(X_i^2) \in\N$. 
\end{proof}
\begin{ex}\emph{
Let $\cat$ be a Tambara-Yamagami (TY) category (see example \ref{ty}). Then $d_+(g)=1$ for $g\in G$. 
Also, $X^2=\sum_{g\in G}g$, so $d_+(X)=\sqrt{\vert G\vert}$. Thus $d_+(\cat)=2\vert G\vert$. \\
In the particular case of the Ising model ($G=\Z_2$), $d_+(\bold 1)=d_+(g)=1$ and $d_+(X)=\sqrt2$, 
and  $d_+(\cat)=4$. 
}\end{ex}

\subsection{Global and FP dimensions}

Until the end of the paper, and without precision, we will assume that our categories are over $\C$. 

\subsubsection{Comparison of global and FP dimension}

Let $\cat$ be a fusion category. 
\begin{thm}[\cite{ENO}]
For every simple object $V$ in $\cat$, one has $\vert V\vert^2\leq d_+(V)^2$, and hence 
${\rm dim}\cat\leq d_+(\cat)$. Moreover, if ${\rm dim}\cat=d_+(\cat)$, then 
$\vert V\vert^2=d_+(V)^2$ for any simple $V$. 
\end{thm}
\begin{proof}
It is sufficient to consider the pivotal case (otherwise one can take the pivotal extension 
$\overline{\cat}$ and recall that the forgetful functor $F:\overline{\cat}\to\cat$ preserves squared 
norms and FP dimension, because it is tensor). 

In this case $N_i\vec d=d_i\vec d$ (where $d_i=\textrm{dim}X_i$ and $\vec d=(d_0,\dots,d_{n-1})$), 
thus by the FP theorem $\vert d_i\vert\leq d_+(X_i)$, and this is an equality if 
$\sum_i\vert d_i\vert^2=\sum_id_+(X_i)^2$. 
\end{proof}
\noindent\textbf{Remark. }
In general, the FP dimension of a fusion category and its global dimension are not equal, or even 
Galois-conjugate (and the same is true for $d_+(V)^2$ and $(\textrm{dim}V)^2$, for any simple object 
$V$). 

Now denote respectively by $D$ and $\Delta$ the global and FP dimensions of $\cat$. We already know 
$D/\Delta\leq1$ (previous theorem), moreover we have 
\begin{thm}[\cite{ENO}]\label{ddel}
$D/\Delta$ is an algebraic integer. 
\end{thm}
\begin{proof}
We can assume $\cat$ is spherical. Otherwise one may consider its pivotal extension, which can be 
shown to be spherical (see \cite{ENO}), and whose global and FP dimensions are respectively $2D$ and 
$2\Delta$). 

In this case $Z(\cat)$ is modular, of global and FP dimensions $D^2$ and $\Delta^2$ (respectively). 
Let $s=(s_{ij})_{ij}$ be  its $S$-matrix. It follows from the Verlinde formula that the matrices 
$N_i$ have common eigenvalues $s_{ij}/s_{0j}$, and the corresponding eigenvectors are the columns of 
$s$. Since $s$ is nondegenerate, there exists a unique label $r$ such that $s_{ir}/s_{0r}=d_+(Y_i)$, 
where $Y_i$ are the simple objects of $Z(\cat)$). \\
Then 
$\Delta^2=\sum_id_+(Y_i)^2=\sum_i\frac{s_{ri}}{s_{0r}}\frac{s_{ir}}{s_{0r}}=\delta_{r^*r}/s_{0r}^2$, 
where we used the  symmetry of $s$ and the fact that $s^2=(\delta_{i^*j})_{ij}$. So we find that 
$r=r^*$ and $\Delta^2=1/s_{0r}^2=D^2/(\textrm{dim}X_r)^2$. Consequently 
$D^2/\Delta^2=(\textrm{dim}X_r)^2$, hence $D/\Delta$ is an  algebraic integer. 
\end{proof}
\begin{cor}[\cite{ENO}]
Let $\cat$ be a nondegenerate fusion category over a field ${\bf k}$ of characteristic $p$. Then its FP 
dimension $\Delta$ is not divisible by $p$. 
\end{cor}
\begin{proof}
Assume that $\Delta$ is divisible by $p$. Let $\widetilde\cat$ be the lifting of $\cat$, and 
$\widehat\cat=\widetilde\cat\otimes_{W(\K)}\mathbb K$ where $\mathbb K$ is the algebraic closure of 
the fraction field of $W(\K)$. Then the Theorem \ref{ddel} says that the global dimension $D$ of 
$\widehat\cat$ is divisible by $\Delta$, hence by $p$. So the global dimension of $\cat$ is zero. 
Contradiction ($\cat$ is nondegenerate). 
\end{proof}

\subsubsection{Pseudo-unitary fusion categories}

\begin{dfn}\emph{
A fusion category $\cat$ (over $\C$) is called \emph{pseudo-unitary} if $\textrm{dim}\cat=d_+(\cat)$. 
}\end{dfn}
\noindent\textbf{Remark. }
Unitary categories (those arising from subfactor inclusions, see \cite{GHJ}) all satisfy this 
condition (so the terminology is coherent). 
\begin{prop}[\cite{ENO}]
Any pseudo-unitary fusion category $\cat$ admits a unique spherical structure, in which 
${\rm dim}V=d_+(V)$. 
\end{prop}
\begin{proof}
Let $b:\textrm{id}\to****$ be an isomorphism of tensor functors, and $g:\textrm{id}\to**$ an 
isomorphism of additive functors such that $g^2=b$. Let $f_i=d_+(X_i)$. Define $d_i=tr(g_{X_i})$ and 
$\vec d=(d_0,\dots,d_{n-1})$; then $f_i=|d_i|$ by pseudounitarity. Further, we can define the action 
of $g$ on ${\rm Hom}(X_i\otimes X_j,X_k)$; let $(T_i)_{jk}$ denote the trace of this operator. Then 
$T_i\vec d=d_i\vec d$, and $|(T_i)_{jk}|\le (N_i)_{jk}$. Thus, 
$$f_if_j=|d_id_j|=|\sum (T_i)_{jk}d_k|\le \sum (N_i)_{jk}f_k=f_if_j\,.$$
This means that the inequality in this chain is an equality. In particular 
$(T_i)_{jk}=\pm (N_i)_{jk}$, and the argument of $d_id_j$ equals the argument of $(T_i)_{jk}d_k$ 
whenever $(N_i)_{jk}>0$. This implies that whenever $X_k$ occurs in the tensor product 
$X_i\otimes X_j$, the ratio $d_i^2d_j^2/d_k^2$ is positive. 
Thus, the automorphism of the identity functor $\sigma$ defined by $\sigma|_{X_i}=d_i^2/|d_i|^2$ 
is a tensor automorphism. Let us twist $b$ by this automorphism, i.e., replace $b$ by $b\sigma^{-1}$. 
After this twisting, the new dimensions $d_i$ will be real. Thus, we can assume without loss of 
generality that $d_i$ were real from the beginning. 

It remains to twist the square root $g$ of $b$ by the automorphism of the identity functor $\tau$ 
given by $\tau|_{X_i}=d_i/|d_i|$ (i.e., replace $g$ by $g\tau$). After this twisting, the new $T_i$ 
is $N_i$ and the new $d_k$ is $f_k$. This means that $g$ is a pivotal structure with positive 
dimensions. It is obvious that such a structure is unique. We are done. 
\end{proof}
\begin{thm}[\cite{ENO}]
Any fusion category of integer FP dimension $\Delta$ is pseudo-unitary. In particular it is 
canonically spherical. 
\end{thm}
\begin{proof}
Let $D=D_1,\dots,D_m$ be the algebraic conjugates of $D=\textrm{dim}\cat$. Then consider 
$g_i\in\textrm{Gal}(\overline\Q/\Q)$ such that $g_i(D)=D_i$, and the corresponding categories 
$\cat_i=g_i(\cat)$. We know that $\textrm{dim}\cat_i=D_i$ and $d_+(\cat_i)=\Delta$, so 
$D_i/\Delta\leq1$ is an algebraic integer. Hence $\prod_i(D_i/\Delta)$ is an algebraic integer 
$\leq1$. But it is also a rational number (because $\prod_iD_i,\Delta\in\N$), so it is an integer 
which is necessarily $1$, and therefore $D_i=\Delta$ for all $i$. In particular $D=\Delta$. 
\end{proof}
\begin{cor}[The Larson-Radford theorem, \cite{LR}]\label{lr}
If $H$ is a finite dimensional semisimple Hopf algebra over $\C$ with antipode $S$, then $S^2=1$. 
\end{cor}
\begin{proof}
Let $\cat=\rep H$. On the one hand we know that $d_+(\cat)=dim(H)\in\N$, hence $\cat$ is 
pseudo-unitary. By example \ref{kap}, it means $dim(H)=\textrm{dim}\cat=tr(S^2)$. \\
On the other hand, $S$ is of finite order, so $S^2$ is semisimple and its eigenvalues are roots of 
unity. Consequently $S^2=1$. 
\end{proof}

\subsection{Classification}

A natural classification problem for fusion categories is the following one. 

\begin{pb}
Classify fusion categories over $\C$ of given Frobenius-Perron dimension. 
\end{pb}
The next theorem solves this problem in the case of the Frobenius-Perron dimension being a prime 
number $p$. Namely, it generalizes to the quasi-Hopf algebra case a result of Kac and Zhu on 
semisimple Hopf algebras of prime dimension $p$. 

Let $\cat$ be a fusion category over $\Bbb C$. 

\begin{thm}[\cite{ENO}]
If $d_+(\cat)=p$ is a prime, then $\cat=\cat(\Z_p,\omega)$. In particular, any semisimple quasi-Hopf 
algebra $H$ of prime dimension $p$ is of the form $H={\rm Fun}(\Z_p)$ with associator defined by 
$\omega\in H^{3}({\Z_p},{\C^\times})=\Z_p$. 
\end{thm}
\begin{proof}
$d_+(\cat)=p$ is a prime, then $d_+(Z(\cat))=p^2\in\N$. Hence $Z(\cat)$ has a canonical spherical 
structure in which $d_i:=\textrm{dim}X_i=d_+(X_i)$ for any simple object $X_i$. Moreover, since 
$\cat$ is itself spherical (because it is of integer FP dimension), $Z(\cat)$ is modular and hence 
$p^2/d_i^2$ is an algebraic integer. Thus $d_i=1$ or $\sqrt p$ (as $d_i^2\in\N$). \\
If there exists $i$ such that $d_i=\sqrt p$, then using the forgetful functor $F:Z(\cat)\to\cat$ we 
find a simple object $F(X_i)$ in $\cat$ with FP dimension $\sqrt p$ (it is simple because the 
dimensions of its simple constituents must be square roots of integers). Since $d_+(\cat)=p$, it is 
the only simple object in $\cat$. This is a contradiction (there must be a neutral object).\\
Consequently, all simple objects in $Z(\cat)$, and hence in $\cat$ also (using $F$), have FP 
dimension $1$, i.e. are invertible. But fusion categories all whose simple objects are invertible 
are all of the type $\cat(G,\omega)$. In our case the group $G$ must have order $p$, so the result 
is proved. 
\end{proof}

With quite a bit more work, this theorem can be extended to the case of products of two primes. 

\begin{thm}
If $d_+(\cat)=pq$ for two prime numbers $p\leq q$, then either $p=2$ and $\cat$ is a Tambara-Yamagami 
category attached to the group $\Z_q$, or $\cat$ is Morita equivalent to $\cat(G,\omega)$ with 
$\vert G\vert=pq$. 
\end{thm}
\begin{proof}
The case $p=q$ is done in \cite[Proposition 8.32]{ENO} and the case $p<q$ is treated in \cite{EGO}. 
\end{proof}

\subsection*{Open problems}

In conclusion we formulate two interesting open problems. 

\begin{enumerate}
\item Let us fix $N\in\N$ (and still work over $\C$). 
E. Landau's theorem (1903) says that the number of finite groups which have $\le N$ irreducible 
representations is finite. In the same way, the number of semisimple finite dimensional quasi-Hopf 
algebras which have $\le N$ irreducible representations is finite (see \cite{ENO}). \\
It is natural to ask if the number of fusion categories over $\C$ with $\le N$ simple objects is 
finite. In the case $N=2$ this is shown in \cite{O2}, but the case $N=3$ is already open. \\
\item Does there exists a semisimple Hopf algebra $H$ over $\C$ whose representation category 
$\rep H$ is not  group-theoretical? \\
For quasi-Hopf algebras, it exists (consider e.g. a TY category related to $G=\Z_p\times\Z_p$ with 
the isomorphism $G^\vee\to G$ corresponding to an elliptic quadratic form, see \cite{ENO}). 
\end{enumerate}


\frenchspacing
\begin{thebibliography}{BaKi}

\bibitem[BaKi]{BK} B.~Bakalov, A.~Kirillov Jr., 
{\it Lectures on tensor categories and modular 
functors,} AMS, Providence 2001. 

\bibitem[B]{B} J.~Bichon, The representation 
category of the quantum group of a non-degenerate 
bilinear form, {\tt arXiv:math.QA/0111114}, 
{\it Comm. Alg.} 31 (2003), 4931--4851. 

\bibitem[BW]{BW} E.~Blanchard, A.~Wassermann, 
in preparation. 

\bibitem[BNS]{BNS} G.~Böhm, F.~Nill, 
K.~Szlach\'anyi, Weak Hopf algebras 
I. Integral theory and $C^*$-strucuture, 
{\it J. Algebra} 221 (1999), 357--375. 

\bibitem[Ca]{Ca} J.~W.~Cassels,
{\it Local fields,} London Math. Soc. 
Student Texts 3, Cambridge Univ. Press, 1986. 

\bibitem[CG]{CG} A.~Coste, T.~Gannon, Remarks on 
Galois symmetry in rational conformal field theories,
{\it Phys. Lett. B} 323 (1994), no. 3--4, 316--321. 

\bibitem[De]{De} P.~Deligne, 
Cat\'egories tensorielles, 
{\it Mosc. Math. J.} 2 (2002), no. 2, 227--248. 

\bibitem[dBG]{dBG} J.~de~Boere, J.~Goeree, 
Markov traces and $II_1$ factors in conformal 
field theory, {\it Comm. Math. Phys.} 139 (1991), 
no. 2, 267--304. 

\bibitem[Dav]{Dav} A.~Davydov, 
Twisting of monoidal structures, 
Max Planck Inst. preprint no. 123, 1995. 

\bibitem[Dr]{Dr} V.G.~Drinfeld, 
Quasi-Hopf algebras, {\it Leningrad Math. J.}
1 (1990), no. 6, 1419--1457. 

\bibitem[EG]{EG} P.~Etingof, S.~Gelaki, The 
classification of triangular semisimple and 
cosemisimple Hopf algebras over an 
algebraically closed field, {\it Intern. Math. 
Res. Notices} 5 (2000), 223--234. 

\bibitem[EG2]{EG2} P.~Etingof, S.~Gelaki, 
On finite dimensional semisimple and cosemisimple 
Hopf algebras in positive characteristic, 
{\it Intern. Math. Res. Notices} 16 (1998), 851--864. 

\bibitem[EG3]{EG3} P.~Etingof, S.~Gelaki, The 
classification of finite dimensional triangular 
Hopf algebras over an algebraically closed 
field of characteristic zero, {\it Moscow Math. J.}
3 (2003), no. 1, 37--43. 

\bibitem[EGO]{EGO} P.~Etingof, S.~Gelaki, V.~Ostrik, 
Classification of fusion categories of dimension $pq$, 
{\tt arXiv:Math.QA/0304194}, {\it Intern. Math. 
Res. Notices} 57 (2004), 3041--3056. 

\bibitem[E]{E} P.~Etingof, On Vafa's theorem 
for tensor categories, {\it Math. Res. Lett.} 
9 (2002), 651--657. 

\bibitem[ENO]{ENO} P.~Etingof, D.~Nikshych, V.~Ostrik, 
On fusion categories, {\tt arXiv:math.QA/0203060}, 
{\it Ann. of Math.} 162 (2005), no. 2, 581--642. 

\bibitem[EO1]{EO1} P.~Etingof, V.~Ostrik, 
Finite tensor categories, {\tt arXiv:math.QA/0301027},
{\it Moscow Math. J.} 4 (2004), no. 3, 627--654. 

\bibitem[EO2]{EO2} P.~Etingof, V.~Ostrik, 
Module categories over representations of $SL_q(2)$ 
and graphs, \texttt{arXiv:math.QA/0302130},  
{\it Math. Res. Lett.} 11 (2004), no. 1, 103--114.

\bibitem[GHJ]{GHJ} F.~Goodman, P.~de la Harpe, 
V.~Jones, {\it Coxeter Graphs and Towers of Algebras,} 
MSRI Publications 14, Springer, 1989. 

\bibitem[HI]{HI} R.~Howlett, and M.~Isaacs, 
On groups of central type, {\it Math. Z.} 179 
(1982), 555--569. 

\bibitem[JS]{JS} A.~Joyal, R.~Street, 
Braided tensor categories, {\it Adv. Math.}
102 (1993), 20--78. 

\bibitem[K]{K} C.~Kassel, {\it Quantum groups}, 
Graduate Texts in Mathematics 155, 
Springer-Verlag, New-York 1995. 

\bibitem[KO]{KO} A.~Kirillov, V.~Ostrik, On a
$q$-analogue of the McKay correspondence and 
ADE classification of $\aff$ conformal field 
theories, {\tt arXiv:math.QA/0101219},  
{\it Adv. Math.} 171 (2002), 183--227. 

\bibitem[LR]{LR} R.~Larson, D.~Radford, 
Semisimple cosemisimple Hopf algebras, 
{\it Amer. J. Math.} 110 (1988), 187--195. 

\bibitem[Mo]{M} M.~Movshev, 
Twisting in group algebras of finite groups, 
{\it Functional Anal. Appl.} 27 (1993), no. 4, 
240--244. 

\bibitem[Mac]{ML} S.~MacLane, 
{\it Categories for the working mathematician} 
(second edition), Graduate Texts in Mathematics 5, 
Springer-Verlag, New-York 1998. 

\bibitem[Mu1]{Mu1} M.~Müger, From subfactors 
to categories and topology I. Frobenius algebras 
in and Morita equivalence of tensor categories, 
{\tt arXiv:math.CT/0111204},  
J. Pure Appl. Alg. 180 (2003), 81--157. 

\bibitem[Mu2]{Mu2} M.~Müger, From subfactors 
to categories and topology II. The quantum double 
of tensor categories and 
subfactors, {\tt arXiv:math.CT/0111205},  
J. Pure Appl. Alg. 180 (2003), 159--219. 

\bibitem[N]{N} D.~Nikshych, 
On the structure of weak Hopf algebras, 
\texttt{arXiv:math.QA/0106010},  
Adv. Math. 170 (2002), no. 2, 257--286. 

\bibitem[NTV]{NTV} D.~Nikshych, V.~Turaev, L.~Vainerman, 
Invariants of knots and 3-manifolds from quantum 
groupoids, {\tt arXiv:math.QA/0006078},  
Topology Appl. 127 (2003), 91--123. 

\bibitem[NV]{NV} D.~Nikshych, L.~Vainerman, 
Finite quantum groupoids and their applications,
\texttt{arXiv:math.QA/0006057},  
in {\it New directions in Hopf algebras,} 
211--262, MSRI Publications 43, 
Cambridge Univ. Press, Cambridge 2002. 

\bibitem[NZ]{NZ} W. D. Nichols and M. B. Zoeller, 
A Hopf algebra freeness theorem, 
{\it Amer. J. Math.} 111 (1989), 381--385. 

\bibitem[O1]{O1} V.~Ostrik, Module categories, 
weak Hopf algebras and modular invariants, 
{\tt arXiv:math.QA/0111139},  
Transform. Groups 8 (2003), 177--206. 

\bibitem[O2]{O2} V.~Ostrik, 
Fusion categories of rank 2, 
{\tt arXiv:math.QA/0203255},  
Math. Res. Lett. 10 (2003), 177--183. 

\bibitem[O3]{O3} V.~Ostrik, 
Module categories over the Drinfeld double of a 
finite group, {\tt arXiv:math.QA/0202130},  
{\it Intern. Math. Res. Notices} 
27 (2003), 1507--1520. 

\bibitem[S]{S} P.~Schauenburg, 
Quotients of finite quasi-Hopf algebras, 
{\tt arXiv:math.QA} {\tt/0204337},  
{\it Hopf algebras in noncommutative geomety and
physics} (Brussels, 2002), 281--290, Lecture Notes 
in Pure and Appl. Math. 231, Marcel Dekker, 
New York 2005.  

\bibitem[TY]{TY} D.~Tambara, S.~Yamagami, Tensor 
categories with fusion rules of self-duality for 
finite abelian groups, 
{\it J. Algebra} 209 (1998), no. 2, 692--707. 

\bibitem[Tu]{T} V.~Turaev, 
{\it Quantum invariants of knots and 3-manifolds}, 
W.~de Gruyter, Berlin 1994. 

\bibitem[Y]{Y} D.~Yetter, 
Braided deformations of monoidal categories and 
Vassiliev invariants, \texttt{arXiv:q-alg/9710010}, 
{\it Higher category theory} (Evanston, IL, 1997), 
Contemp. Math. 230 (1998), 117--134, AMS, 
Providence, RI. 
\end{thebibliography}
\end{document}